\documentclass[10pt,english]{article}
\usepackage[letterpaper,bindingoffset=0.2in,%
            left=0.75in,right=0.75in,top=1in,bottom=1in,%
            footskip=.25in]{geometry}
\usepackage{cite}
\usepackage{notoccite}
\makeatletter
\renewcommand\@biblabel[1]{#1.}
\makeatother
\usepackage{amsthm}
\newlength{\bibitemsep}
\newlength{\bibparskip}
\let\oldthebibliography\thebibliography
\renewcommand\thebibliography[1]{%
  \oldthebibliography{#1}%
  \setlength{\parskip}{\bibitemsep}%
  \setlength{\itemsep}{\bibparskip}%
}
\usepackage{epsfig, xcolor, setspace, graphicx}

\usepackage[leftcaption]{sidecap}
\usepackage{lscape,amsmath,nccmath,nicefrac}
\newcommand\myeq{\mathrel{\overset{\makebox[0pt]{\mbox{\normalfont\tiny\sffamily def}}}{=}}}
\usepackage{relsize}
\usepackage{setspace}
 
\usepackage[T1]{fontenc}
\usepackage{footmisc}
\usepackage{hyperref} 
\hypersetup{colorlinks,breaklinks,
	unicode = true,
            urlcolor=[rgb]{0,0,.5},
            linkcolor=[rgb]{0,0,.5},
            citecolor=[rgb]{0,0,.5}}
\usepackage[all]{hypcap} 
\usepackage{caption}
\usepackage{makecell,multirow}
\usepackage{hhline}
\usepackage[normalem]{ulem} 
\usepackage{float}
\usepackage{centernot}
\usepackage{etoolbox}
\usepackage{amssymb}
\usepackage{tabularx}
\usepackage{wrapfig,booktabs}

\sidecaptionvpos{figure}{t}

\apptocmd{\lim}{\limits}{}{}
\newcolumntype{C}[1]{>{\centering\arraybackslash}p{#1}}

\newcommand{\cmmnt}[1]{\ignorespaces}
\newcolumntype{Y}{>{\centering\arraybackslash}X} 
\usepackage[euler]{textgreek}
\usepackage{cleveref}

\usepackage{tikz}
\usetikzlibrary{shapes.geometric, arrows, shapes, shapes.multipart}

\title{A series acceleration algorithm for the gamma-Pareto (type I) convolution and related functions of interest for pharmacokinetics \vspace{1.5cm}}

\author{{\bf Carl A. Wesolowski$^{\dag\ddag\footnote{Corresponding Author: telephone :(306) 665 1515, e-mail: carl.wesolowski@gmail.com}}$, Jane Alcorn$^{\dag}$, and Geoffrey T. Tucker$^{\S}$} \vspace{0.5cm}\\ 
       $^{\dag}$  {\normalsize College of Pharmacy and Nutrition}\\
       {\normalsize University of Saskatchewan, 104 Clinic Place, Saskatoon, SK, S7N 2Z4, Canada}\\
       $^{\ddag}$ {\normalsize Department of Medical Imaging, Royal University Hospital, College of Medicine}\\
       {\normalsize University of Saskatchewan, 103 Hospital Drive, Saskatoon, SK, S7N 0W8, Canada}\\
       $^{\S}$ {\normalsize Department of Human Metabolism, Medical and Biological Sciences, University of Sheffield, Sheffield, UK}}
    
\date{}
\begin{document}  

\tikzstyle{beginend} = [rounded rectangle, minimum width=2.0cm, minimum height=1cm,text centered, draw=black ]
\tikzstyle{decision} = [diamond, aspect=2, minimum width=2.5cm, minimum height=.8cm, text centered, text height=0.25cm, draw=black ]
\tikzstyle{io} = [trapezium, trapezium left angle=70, trapezium right angle=110, minimum width=2.5cm, minimum height=0.8cm, text centered, text width=1.5cm, draw=black ]
\tikzstyle{process} = [rectangle, minimum width=2.cm, minimum height=1cm, text centered, text width=2.0cm, text height=0.25cm, draw=black ]
\tikzstyle{arrow} = [thick,->,>=stealth ]
\tikzstyle{pp} = [rectangle split, rectangle split horizontal, rectangle split parts=3, minimum height=1cm, text centered, text width=2cm, text height=0.4cm ]

\begin{onehalfspacing}
\date{Accepted for Publication, Aug. 26, 2021 https://doi.org/10.1007/s10928-021-09779-4} 
\maketitle \vspace{0in} \noindent




\begin{abstract}
\noindent The gamma-Pareto type I convolution (GPC type I) distribution, which has a power function tail, was recently shown to describe the disposition kinetics of metformin in dogs precisely and better than sums of exponentials. However, this had very long run times and lost precision for its functional values at long times following intravenous injection. An accelerated algorithm and its computer code is now presented comprising two separate routines for short and long times and which, when applied to the dog data, completes in approximately 3 minutes per case. The new algorithm is a more practical research tool. Potential pharmacokinetic applications are discussed.

 \vspace{2em}

\noindent $\mathbf{Keywords}$: Algorithm, Gamma-Pareto convolution, Metformin, Pharmacokinetics, Statistical distributions; Mathematica
\end{abstract}

\def\thesection{\Roman{section}}
\thispagestyle{plain}

\section*{Introduction} Gamma-Pareto convolutions (GPC), the convolution of a gamma distribution with some type of Pareto distribution, are increasingly used for modelling diverse random processes like traffic patterns \cite{nadarajah2007convolution}, flood rates, fatigue life of aluminium, confused flour beetle populations \cite{alzaatreh2012gamma}, and extreme rainfall events \cite{hanum2015modeling}. 
Although there are multiple possible GPC models and different nomenclatures used to describe them, a natural classification would arise from Pareto distribution classification, types I through IV, and the Lomax distribution, a type II subtype, which is the classification scheme of reference \cite{kotz2004continuous} and the Mathematica computer language \cite{Mathematica}.\footnote{Wolfram Research, Inc., (2021) Mathematica, Version 12.3, Champaign, IL \\\indent https://reference.wolfram.com/language/ref/ParetoDistribution.html} 

Convolution was first introduced to pharmacokinetics in 1933 by Gehlen who used the convolution of two exponential distributions,

\begin{equation}
\text{EDC}(b,\beta ;t)=b\, e^{-b\,x}* \beta e^{-\beta\,x}\,(t)\\
=
 \begin{cases} 
 \begin{array}{ll} b\, \beta \frac{ e^{-\beta t}-e^{-b\, t}} {b-\beta }&b\neq \beta\\
b^2 t\, e^{-b\, t}&b=\beta
\end{array}
\Bigg\}
 & t\geq 0\\
\;\;0\hspace{9em}\}&t<0
\end{cases}\;,
\end{equation}

\noindent to describe plasma concentration-time data, and as originally developed in 1910 by Bateman to model radioactive decay \cite{bateman1910solution,gladtke1988history,gehlen1933wirkungsstarke}. 
Much later, in 2006, the Bateman equation was generalised as an exact gamma-gamma convolution (GDC) by Di Salvo \cite{di2006exact}. Ten years later, this was then applied to 90 min continuous recordings of radioactivity in human thyroid glands following injection of $^{99m}$Tc-MIBI \cite{Wesolowski2016GDC}.\footnote{740 MBq technetium-99m labeled hexakis-methoxy-isobutyl-isonitrile.} In 1919, Widmark identified integration of a monoexponential as a model for constant infusion \cite{widmark1919studies}. That integration from zero to $t$ to find a constant infusion model can be applied not just to exponential functions, but applies equally well to any area under the curve (\textit{AUC}) scaled density-function (pdf)\footnote{\label{note1}We retain the acronym pdf without a probability;\textit{ p}, but use $f(t)$ preferentially. Concentration models are the product of area-under-the-curve of concentration and density functions whose total area-under-the-curve is 1 (dimensionless dose fraction). This balances the classical mechanical units of Mass, Length, and Time, as follows,
$C(t)=\textit{AUC}\, \times\,f(t) \;;\;\;\;\left[ \frac{M}{L^{3}}=\frac{M\;T}{L^{3}}\;\times \,\frac{1}{T}\right]\,.$} model, was shown subsequently by Wesolowski \textit{et al.} \cite{Wesolowski2016PLoS,Wesolowski_2020}. Recently, the disposition of metformin was described precisely using the type I GPC model, which as it is asymptotically a Pareto function, has a power function tail \cite{Wesolowski_2020}. Using direct comparison rather than classification, power function tails were shown to be always heavier than exponential tails, see the Appendix Subsection entitled \textit{Relative tail heaviness} of reference \cite{Wesolowski_2020}. A power function tail, in turn, implies an underlying fractal structure, where fractal in this context signifies a scale invariant model of vascular arborisation \cite{west1999fourth}. 

The GPC computer algorithm used in 2019 had long run times and was not accurate beyond 96 h \cite{Wesolowski_2020}. These problems were corrected in order to make predictions for multiple dosing over longer times \cite{Tucker2020,Tucker2020a}. Since computer implementation of new functions is highly specialised, not easily arrived at by induction and yet indispensable for any practical application, documentation of a more practical type I GPC algorithm may facilitate its more widespread implementation. Accordingly, we now present a series acceleration computer implementation of a more generally applicable GPC type I function, with markedly shorter run times.

\section*{Background}

\subsection*{The gamma-Pareto convolution distribution function family}

A classification for gamma-Pareto convolutions (GPC) is proposed that arises from the types of Pareto distributions \cite{kotz2004continuous}. These are types I through IV plus the Lomax distribution, a subtype of II. The Pareto type I distribution is
 \begin{equation}\label{eq:PD}
\textnormal{PD}(t; \alpha, \beta)=
 \dfrac{\alpha}{t} \left(\dfrac{\beta}{t}\right) ^{\alpha } \theta(t-\beta)\;,
 \end{equation} 
 \noindent where $\alpha$ is the shape parameter, $\beta$ is a scale parameter and $\theta(\cdot)$ is the unit step function such that $\theta(t-\beta)$ is the unit step function time-delayed by $\beta$, and is used to make a product that is non-zero only when $t> \beta$.\footnote{The unit step function, $\theta(x)$, is zero for $x<0$ and 1 for $x\geq 0$, such that $\theta(x)$ is continuous everywhere except at $x= 0$. When $x=t-\beta$ and $\beta>0$, then $\theta(t-\beta)$ is a unit step function shifted to later time (i.e., to the right) by $\beta$ units in the new coordinate system; $t$. The unit step function is faster for numerical computations than the Heaviside theta function, which later is sometimes also symbolised as $\theta(x)$. The Heaviside theta is more mathematically useful when it is continuous everywhere such that its derivative and Laplace transform are defined.} 

A type II Pareto distribution can be written as 

\begin{equation}\text{PD}_{\text{II}}(t;\alpha,\beta,\mu)=\frac{\alpha }{\beta }\left(1+\frac{t-\mu }{\beta }\right)^{-\alpha -1}\theta(t-\mu)\;.\end{equation}

\noindent Setting $\mu=0$, this becomes the Lomax distribution; $\text{PD}_{\text{Lomax}}(t;\alpha,\beta)=\frac{\alpha }{\beta }\big(1+\frac{t }{\beta }\big)^{-\alpha -1}\theta(t),$ which was used to derive a Lomax gamma-Pareto distribution \cite{nadarajah2007convolution}. The relevance of this is that the GPC type I and Lomax GPC derivations are similar. As yet, the type II (not Lomax) through type IV gamma-Pareto convolutions have not been published. These convolutions are likely to be infinite sums and may require series acceleration to be of practical use. By substitution and reduction of the number of parameters, there are closed form GPC-like expressions, types II through IV, that are different distributions \cite{alzaatreh2012gamma}. As a full set of solutions for the entire GPC function family has not been characterised, it is not known what additional applications there could be for the GPC family of functions. 

Unlike the Lomax GPC, the GPC type I does not start at $t=0$, but at $t=\beta$. For pharmacokinetic modelling, $\beta>0$ is a measure of the circulation time between injection of an intravenous bolus of drug ($t=0$), and its arrival at a peripheral venous sampling site ($t=\beta$). The four-parameter gamma Pareto (type I) convolution (GPC) density function was developed to model the disposition of metformin in dogs, which exhibited an unexpectedly heavy tail poorly described by an exponential decay \cite{Wesolowski_2020}. This heavy tail implies a prolonged buildup of the body burden of the drug \cite{tucker1981metformin,Tucker2020} that may require dose tapering on long-term use, especially in patients with renal impairment \cite{Tucker2020a}.

\subsection*{The Gamma-Pareto type I convolution and related functions}

\textbf{GPC type I:} To form a GPC type I model, the type I Pareto distribution, Eq.~\eqref{eq:PD}, is convolved with a gamma distribution,
\begin{equation}\label{eq:GD}
 \text{GD}(t; a,b) =
 \,\dfrac{1}{t}\;\dfrac{e^{-b \, t}(b \, t)^{\,a} }{\Gamma (a)}\theta(t)\;,
 \end{equation}

\noindent where $a$ is a dimensionless shape parameter, $b$ is a rate per unit time, is the reciprocal of its scale parameter, and $\Gamma(\cdot)$ is the gamma function.\footnote{The gamma function, or generalised factorial, is $\Gamma(z)=\int_0^{\infty } \frac{t^{z-1}}{e^t} \, dt;\, \Re(z) >0$} This yields the GPC function,

\begin{equation}\label{eq:GPC}
\text{GPC}(t)=\theta (t-\beta)\;\frac{b^a\, \alpha\, \beta ^{\alpha } }{\Gamma(a)}t^{a-\alpha -1}\sum _{n=0}^{\infty } \frac{(-b\, t)^n }{n!}B_{1-\frac{\beta }{t}}\left(a+n,-\alpha \right)\;,
\end{equation}

\noindent where $B_z(\cdot,\cdot)$ is the incomplete beta function.\footnote{The incomplete beta function is $B_z(a,b)=\int_0^z t^{a-1} (1-t)^{b-1} \, dt;\,\Re(a)>0\land \Re(b)>0\land | z| <1$} This is a density function (a pdf, or more simply an $f$, with units per time; $t^{-1}$). Equation \eqref{eq:GPC} is from convolution following Maclaurin series expansion of $e^{-b\,t}$, i.e., it is analytic. An analytic function has any number of sequential multiple integrals and derivatives, as illustrated in the following equations. Compared to their prior expression \cite{Wesolowski_2020}, the equations that follow have been put in simpler terms. 
\newline
\newline
\noindent\textbf{GPC type I integral:} Equation \eqref{eq:CDF} is the cumulative density function (CDF) of the GPC, symbolised by $F$, the integral of the $f(t)$ density; $F(t)=\int_0^t f(\tau ) \, d\tau$,

\begin{equation}\label{eq:CDF}
\text{GPC}_{F}(t) =\theta (t-\beta)\;\frac{b^a \alpha\, \beta ^{\alpha }}{\Gamma(a)}t^{a-\alpha } \sum _{n=0}^{\infty} \frac{(-b\, t)^n}{(a+n) n!} B_{1-\frac{\beta }{t}}\left(1+a+n,-\alpha\right)\;.
\end{equation}

\noindent This equation, because it is a CDF, expresses the dimensionless fraction of a unit drug dose eliminated from the body as a function of time, and was used to calculate a prolonged retention of metformin in dogs and to explain its incomplete urinary recovery at 72 h following intravenous injection in humans \cite{Wesolowski_2020,tucker1981metformin,Tucker2020}. 
\newline
\newline
\textbf{GPC type I double integral:} Equation \eqref{eq:SCD} is the double integral of the density function, $f$, which is also the single integral of $F$, the CDF, and is sometimes called a "super-cumulative" distribution \cite{Avdis_2017}. It is symbolised by $\mathcal{F}$, i.e., $\mathcal{F}(t)=\int_0^t F(\tau )\, d\tau=\int_0^t \int_0^\tau f(x) \,d x\, d\tau $. The GPC$_{\mathcal{F}}$ in least terms is 

\begin{equation}\label{eq:SCD}
\text{GPC}_{\mathcal{F}}(t)=\theta (t-\beta)\;\frac{b^a \alpha\, \beta ^{\alpha }}{\Gamma(a)}t^{a-\alpha +1} \sum _{n=0}^{\infty } \frac{(-b\, t)^n\text{  }}{(a+n)(1+a+n) n!}B_{1-\frac{\beta }{t}}\left(2+a+n,-\alpha \right)\;.
\end{equation}

\noindent This equation (units $t$) was used to construct an intravenous bolus multidose loading regimen that maintains the same mean amount of metformin in the body during successive dose intervals \cite{Wesolowski_2020} and to predict metformin buildup during constant multidosing in humans both with normal renal function and with renal insufficiency \cite{Tucker2020a}. A further use of this equation is to predict the cumulative distribution function following a period of constant infusion given only its bolus intravenous-concentration, fit function.
\newline
\newline
\noindent\textbf{GPC type I derivative:} Equation \eqref{eq:dgpc} is the derivative of the GPC density, GPC$'$, or in general an $f'$,

\begin{equation}\label{eq:dgpc}
\text{GPC}'(t)=\theta (t-\beta)\;\frac{b^a \alpha\, \beta ^{\alpha } }{\Gamma(a)}t^{a-\alpha -2}\sum _{n=0}^{\infty } \frac{(a+n-1) (-b \,t)^n }{n!}B_{1-\frac{\beta}{t}}\left(a+n-1,-\alpha \right)\;.
\end{equation}

\noindent This equation (units $t^{-2}$) is useful for finding the peaks of the GPC function by searching for when it equals zero, and for calculating disposition half-life from its general definition, 
\[t_{1/2};f(t) \myeq -\ln(2)\mfrac{f(t)}{f'(t)}\;,\]
which is Eq.~(6) of reference \cite{Wesolowski_2020}. Note that there is a pattern in the sequential integrals and derivatives that illustrates the analyticity of the GPC function. The integrals and derivatives above follow directly from integration or differentiation of the GPC formula, for which the following identity from integration by parts\footnote{The parts are $U(x)=\frac{1}{\text{A}}\big(\frac{1}{1-x}\big)^{-\text{A}} (1-x)^\text{B}$, and $V(x)=\big(\frac{x}{1-x}\big)^\text{A}$. The identity is listed elsewhere: Wolfram Research Inc. (2021), Champaign, IL. http://functions.wolfram.com/06.19.17.0001.01.} $$B_z(\text{A}+1,\text{B})=\frac{\text{A}}{\text{A}+\text{B}}B_z(\text{A},\text{B})-\frac{z^\text{A} (1-z)^\text{B}}{\text{A}+\text{B}}\;\;,$$ is useful for simplifying the results.

\section*{Methods, algorithms for GPC type I series acceleration and their computation}

\subsection*{Data sources and regression methods}\label{data}
The source data for regression analysis and subsequent algorithm testing consists of seven intravenous bolus metformin studies performed in healthy mixed-breed dogs \cite{johnston2017pharmacokinetics}. The 19 to 22 samples per case drawn between 20 min to 72 h postinjection are listed as open data in \textit{Supplementary material 1 (SLSX 49kb)} in \cite{Wesolowski_2020}.\footnote{https://link.springer.com/article/10.1007/s10928-019-09666-z\#Sec220} The regression target was the so-called $1/C^2$ weighted ordinary least squares (OLS) method, implemented as minimisation of the proportional norm, which is also the relative root mean square (rrms) error, as per the \textit{Concentration data and fitting it} Appendix Subsection of \cite{Wesolowski_2020}.\footnote{https://link.springer.com/article/10.1007/s10928-019-09666-z\#appendices} The loss function chosen to be minimised agreed with the error type the measurement system assay calibration curve. Both the metformin assay (5.2\% rrms) \cite{michel2015}, and the GPC residuals (8.6\% rrms) exhibited proportional error. The reuse of assay loss functions for regression loss functions is systemically consistent and appears in these references \cite{Wesolowski_2020, Wesolowski2016GDC}. The regression method used was Nelder-Mead \textit{Constrained Global Numerical Minimisation} as implemented in Mathematica, a global search technique \cite{Mathematica}.\footnote{https://reference.wolfram.com/language/tutorial/ConstrainedOptimizationGlobalNumerical.html\#252245038} For 20 significant figure results for all parameters used was the Mathematica routine NMinimize with the options: PrecisionGoal $\to$ 30, \mbox{AccuracyGoal $\to$ 30,} \mbox{WorkingPrecision $\to$ 65,} \mbox{MaxIterations $\to$ 20010,} Method $\to$ \{"NelderMead", "PostProcess" $\to$ False\}. Post processing is disallowed because it launches a constrained convex gradient solution refinement protocol; the interior point method, which does not converge. The use of parameter starting value ranges close to the solution helps speed up convergence. Note that regression can start with 65 significant figure accuracy but finish with less than half of that for some parameter values due to error propagation from the fit function itself and/or the regression process. In order to calculate the confidence intervals (CI) of the parameters, model-based bootstrap \cite{bollen1992bootstrapping} was performed, as follows. Care was taken to verify the normality of fit residuals and the homoscedasticity of residuals---see \cite{Wesolowski_2020}---as suggested by \cite{zhang2016bootstrapping}. Those conditions allow for the residuals to be randomly sampled with replacement, then added to the model at the sample-times to create synthetic data having the same properties as the original data, but which have altered regression parameter solutions. The bootstrap parameter values so obtained can provide more information than gradient method parameter CV's, as the latter only provides single case-wise estimates, which are not as statistically useful as case-wise distributed parameter information \cite{Wesolowski_2020}. Table \ref{params} shows both case-wise and population-wise coefficients of variation from an early version of a GPC algorithm. The table was amalgamated from Tables 1, 3, and 12 of \cite{Wesolowski_2020} representing 24 h of 8-core parallel processing of 42 time-sample serum curves. There is thus an obvious need for a faster algorithm for regression analysis.

 \begin{table}[ht]
\centering
\captionsetup{justification=justified,margin=0cm}
 \caption {Shown are parameters from gamma densities (GD), Pareto densities (PD) and both from Gamma-Pareto convolution (GPC) fitting of concentrations data for 7 dogs with model-based bootstrap root mean square case-wise (Case\%, $n=5$) and population-wise (Pop.\%, $n=35$) coefficients of variation, and fit error.}
 \label{params} 
 \begin{tabularx}{\textwidth}{lccccccccccccc}
 \hline
\small{Functions}\!\!&\multicolumn{4}{c}{GD} &\multicolumn{4}{c}{PD} &\multicolumn{5}{c}{GPC}\\
\hhline{~----~~~~-----} 
\small{Parameters}\!\!& $a$ && $b$ && $\alpha $ && $\beta $ && \emph{AUC} && $\mathit{CL}$ &&Fit error\\
\hhline{-~~~~----~~~~~}
Units$^\text{ a}$ & none &\%&$\frac{1}{\text{h}}$&\%& none &\%&s &\%& $\frac{\text{mg}\cdot \text{h}}{\text{L}}$ &\%& \hspace*{-.4em} $\frac{\text{ml}}{\text{min}\cdot \text{kg}}$ &\% &\%\\[.2em]
\Xhline{2\arrayrulewidth} 
Dog 1 & 0.3493 &8.09& 0.7318 &4.84& 0.2644 &5.13& 25 &---& 31.16 &6.22&  9.8 &6.42 &8.7\\
Dog 2 & 0.8112 &10.5& 0.9993 &6.89& 0.1365 &4.90& 25 &---& 28.18 &1.09& 11.5 &1.09&6.3 \\
Dog 3 & 0.6689 &8.67& 0.9107 &5.95& 0.2010 &3.81& 25 &---& 12.15 &3.77& 26.7 &3.72&5.9\\
Dog 4 & 0.6092 &22.3& 0.8062 &15.9& 0.1726 &9.80& 25 &---& 16.73 &5.65& 19.4 &5.89&13.8\\
Dog 5 & 0.6435 &20.6& 1.1035 &11.4& 0.1199 &6.11& 25 &---& 26.21 &5.37& 12.4 &5.14&9.5\\
Dog 6 & 0.5194 &7.42& 0.6137 &4.93& 0.1929 &5.85& 30 &---& 28.43 &2.15& 11.4 &2.10&6.1\\
Dog 7 & 0.7629 &17.7& 1.0518 &20.3& 0.1571 &5.82& 30 &---& 22.10 &2.31& 14.7 &2.34&10.0\\

Case\% & &14.8&  &11.5&  &6.17&  & --- & &4.22&  &4.26&8.2\\
Pop.\% && 29.7 && 25.0 && 25.9 && --- && 28.8 &&  39.5 &7.5$^{\text{ b}}$\\

\hline
\end{tabularx}
 \begin{tabularx}{1\textwidth}{X} 
 $^\text{a }$Units row: \emph{None} means \emph{dimensionless}. As $\beta$ was constrained to be within 25 to 30 s, its variability for the 5 realisations per case is not meaningful.\\
 $^\text{b }$Geometric mean (GM) was used to calculate group error of fit because the 35 model-based bootstrap fit errors errors were log-normally distributed for which the GM was 7.5\%, not significantly different from the original data GM error of 8.2\%. The original data values are case-wise listed for fit error, and for the parameter values.

\end{tabularx}
\end{table}

\subsection*{GPC type I primary definition: The short-\textit{t} algorithm}
\noindent The primary definition of a gamma-Pareto type I convolution, Eq.~\eqref{eq:GPC}, is

\begin{equation}\label{eq:GPC2}
\begin{split}
\textnormal{GPC}\arraycolsep=1.2pt\def\arraystretch{.7}
\left(\begin{array}{cc}
a&b\\
\alpha&\beta
\end{array} \Big|\,t\right)&=\textnormal{GD}( a,b;x)\ast \textnormal{PD}(\alpha , \beta;x) \;(t)\\
&=\left[\dfrac{(b \, x)^{\,a}\,e^{-b \, x} }{x\,\Gamma (a)}\theta(x)\right]* \left[\dfrac{\alpha}{x} \left(\dfrac{\beta}{x}\right) ^{\alpha } \theta(x-\beta)\right]\;(t)\\
&=\theta (t-\beta)\;
 t^{a-\alpha-1} \frac{\alpha\, b^a\beta^\alpha}{\Gamma (a)}\sum _{n=0}^{\infty } \frac{(-b\,t)^n}{n!} B_{1-\frac{\beta}{t}}(a+n,-\alpha)
\end{split}\;\;. 
\end{equation}

\noindent This contains alternating terms in the summation such that the sum is rapidly convergent for $t$ not much greater than its lower limit, $\beta$. However, for sufficiently large values of $t$, the individual terms of the summation both alternate in sign and become extremely large in \textit{magnitude} (i.e., absolute value) before absolute series convergence. For absolute convergence of an alternating series the infinite sum of the absolute values is bounded above, which permits rewrite of the summation sequence of infinite sums. This, and the ratio test \cite{laugwitz1994riemann} for it are shown in the \nameref{short} Appendix Subsection. Thus, the order of infinite summation can be changed to obtain shorter run times when $t\gg\beta$, and the algorithm is accelerated through an algebraic rewrite of Eq.~\eqref{eq:GPC2} as Eq.~\eqref{eq:accelgpc} below. Alternating infinite series with large magnitude terms occurring before absolute convergence are common, for example, the infinite-series, primary definition of $\sin(x)=x-\frac{x^3}{3!}+\frac{x^5}{5!}-\frac{x^7}{7!}+\frac{x^9}{9!}-\cdots$ has that same property for larger magnitude $x$-values. Acceleration for the sine function could include converting the $x$-values to be principal sine values ($-\frac{\pi}{2}$ to $\frac{\pi}{2}$), and adjusting the output accordingly.\footnote{$\sin(12)$ executes to 65 decimal places in 19 microseconds in the Mathematica language on an 2.3 GHz 8-Core Intel Core i9 processor. Current acceleration algorithms for routine functions are many generations beyond what is outlined here.} For the GPC$(t)$ function a similar result, i.e., adjusting the algorithmic behaviour to be accelerated for long-$t$ values, can be obtained as follows.

\subsection*{GPC type I secondary definition: The long-\textit{t} algorithm}

\newtheorem{1em}{Theorem}
\begin{1em}\label{longT}The long-$t$ algorithm is
 \begin{equation}\label{eq:accelgpc}
\begin{split}
\textnormal{GPC}\arraycolsep=1.2pt\def\arraystretch{.7}
\left(\begin{array}{cc}
a&b\\
\alpha&\beta
\end{array} \Big|\,t\right)=&-\theta (t-\beta) \frac{ \alpha b^a }{\Gamma (a)}t^{a-1}\sum _{k=1}^{\infty }\left(\mfrac{\beta }{t}\right)^k \frac{(1-a)_k}{k! (k-\alpha )} \, _1F_1(a,a-k;-b t)\\
&+\theta (t-\beta )\left[ \frac{b^a }{\Gamma (a)} e^{-b t} t^{a-1}-\pi \csc (\pi \alpha )\frac{b^a \beta ^{\alpha } }{\Gamma (\alpha )} t^{a-\alpha -1} \, _1\tilde{F}_1(a,a-\alpha ;-b t)\right]
\end{split}\;\;\;.
\end{equation}\end{1em}

\noindent \begin{proof} This is shown by substitution of the identities,\footnote{$\,$http://functions.wolfram.com/06.19.17.0008.01 and http://functions.wolfram.com/06.18.02.0001.01} $B_z(A,B)=B(A,B)-B_{1-z}(B,A)$, and $B(A,B)=\frac{\Gamma (A) \Gamma (B)}{\Gamma (A+B)}$ into the incomplete beta function of Eq.~\eqref{eq:GPC2} above and yields,

\begin{equation}\label{eq:3}
B_{1-\frac{\beta }{t}}(a+n,-\alpha )=\frac{\Gamma (-\alpha ) \Gamma (a+n)}{\Gamma (a+n-\alpha )}-B_{\frac{\beta }{t}}(-\alpha ,a+n)\;.
\end{equation}

\noindent Substituting this into the right hand side of Eq.~\eqref{eq:GPC2} yields,

\begin{equation}\label{eq:4}
\theta (t-\beta ) \frac{\alpha\, b^a \beta ^{\alpha } }{\Gamma (a)}t^{a-\alpha -1} \sum _{n=0}^{\infty } \left[\frac{(-b t)^n} {n!}\frac{(\Gamma (-\alpha ) \Gamma (a+n))} {\Gamma (a+n-\alpha )}-\frac{(-b\, t)^n}{n!}B_{\frac{\beta }{t}}(-\alpha ,a+n)\right]\;\;,
\end{equation}

\noindent the left hand summand of which simplifies to a GPC asymptote for long times, $t$,

$$
\theta (t-\beta )\alpha\, b^a \beta^\alpha \Gamma (-\alpha) t^{a-\alpha-1} \, _1\tilde{F}_1(a,a-\alpha;-b t)\;,
$$

 \noindent where $ \, _1\tilde{F}_1(\cdot,\cdot ;z)$ is the regularised confluent hypergeometric function.\footnote{where $\, _1\tilde{F}_1(a;b;z)=\, _1F_1(a;b;z) /\Gamma (b)$, where $\, _1F_1(a;b;z)=\sum _{k=0}^{\infty } z^k (a)_k/[k! (b)_k]$ is the not regularised version, and where $ (a)_k=\Gamma (a+k)/\Gamma (a)$ is the Pochhammer, also called the descending factorial.} The above formula, as\footnote{http://functions.wolfram.com/06.05.16.0001.01 and http://functions.wolfram.com/06.05.16.0002.01} $-\pi \csc (\pi \alpha )=\Gamma (-\alpha ) \Gamma (\alpha +1)=\alpha \Gamma (-\alpha ) \Gamma (\alpha )$, can be written alternatively as
 
\begin{equation}\label{eq:asy}-\theta (t-\beta ) \pi \csc (\pi\, \alpha )\frac{b^a \beta ^{\alpha } }{\Gamma (\alpha )} t^{a-\alpha -1} \, _1\tilde{F}_1(a,a-\alpha ;-b t)\;,
\end{equation}

\noindent which obviates having to use a $\Re[\Gamma(-\alpha)]$ computer command to truncate a zero magnitude imaginary machine number carry, e.g., $\Re(x+0\times i)=x$, such that Eq.~\eqref{eq:GPC2} can be rewritten as 

\begin{equation}\label{eq:proved}
\begin{split}
\textnormal{GPC}\arraycolsep=1.2pt\def\arraystretch{.7}
\left(\begin{array}{cc}
a&b\\
\alpha&\beta
\end{array} \Big|\,t\right)
=&-\theta (t-\beta) \frac{\alpha\, b^a \beta^\alpha}{\Gamma (a)}t^{a-\alpha-1} \sum _{n=0}^{\infty } \frac{(- b\,t)^{n} } {n!} B_{\frac{\beta}{t}}(-\alpha,a+n)\\
&-\theta (t-\beta ) \pi \csc (\pi\, \alpha )\frac{b^a \beta ^{\alpha } }{\Gamma (\alpha )} t^{a-\alpha -1} \, _1\tilde{F}_1(a,a-\alpha ;-b\, t)
\end{split}\;\;\;,
\end{equation}

\noindent where $\alpha\neq0,1,2,3,\dots$, which is Eq.~(25) of the first type I GPC publication \cite{Wesolowski_2020}. Note that not only is the summation of the above absolutely convergent, but as the second line above is an asymptote for $t\to \infty$ of the GPC function \cite{Wesolowski_2020}, the summation converges to zero as $t\to \infty$ relatively more rapidly than the asymptote.

The summation terms, 
$$-\frac{\alpha\, b^a \beta^\alpha}{\Gamma (a)}t^{a-\alpha-1} \sum _{n=0}^{\infty } \frac{(- b\,t)^{n} } {n!} B_{\frac{\beta}{t}}(-\alpha,a+n)\;,$$

\noindent are rearranged for acceleration at long times using the infinite series definition of the incomplete beta function,\footnote{$\,$http://functions.wolfram.com/06.19.06.0002.01}

\begin{equation}\label{eq:betaidentity}
B_{\frac{\beta }{t}}(-\alpha ,a+n)=\left(\frac{\beta }{t}\right)^{-\alpha } \sum _{k=0}^{\infty } \frac{\left(\frac{\beta }{t}\right)^k (1-a-n)_k}{k! (k-\alpha )}\text{ for }\left| \frac{\beta }{t}\right| <1\text{ and } \alpha\neq0,1,2,3,\dots\;,
\end{equation}

\noindent by substituting it into the summation, and simplifying to yield,

\begin{equation}\label{eq:re1}
-\frac{\alpha\, b^a }{\Gamma (a)}t^{a -1}\sum _{n=0}^{\infty } \frac{(-b\, t)^n }{n!}\sum _{k=0}^{\infty } \frac{\left(\frac{\beta }{t}\right)^k (1-a-n)_k}{k! (k-\alpha )}\;.
\end{equation}

\noindent Given absolute convergence (\nameref{short} Appendix Subsection) the order of infinite summation can be changed with impunity by distributing the outer sum over the inner sum, and factoring, as follows,
 
\begin{equation*}-\frac{\alpha\, b^a }{\Gamma (a)}t^{a -1}\sum _{k=0}^{\infty } \sum _{n=0}^{\infty } \frac{(-b\, t)^n }{n!}\frac{\left(\frac{\beta }{t}\right)^k (1-a-n)_k}{k! (k-\alpha )}\;,\end{equation*}
\begin{equation}\label{eq:re2}-\frac{ \alpha\, b^a }{\Gamma (a)}t^{a-1}\sum _{k=0}^{\infty } \frac{\left(\frac{\beta }{t}\right)^k}{k! (k-\alpha )}\sum _{n=0}^{\infty } \frac{(-b\, t)^n (1-a-n)_k}{n!}\;.\end{equation}

\noindent Fortunately, the inner sum in the above formula simplifies to a closed form, allowing it to be rewritten as
\begin{equation}\label{eq:few}
-\frac{ \alpha\, b^a }{\Gamma (a)}t^{a-1}\sum _{k=0}^{\infty } \frac{\left(\frac{\beta }{t}\right)^k}{k! (k-\alpha )}(1-a)_k \, _1F_1(a,a-k;-b\, t)\;.
\end{equation}

\noindent The $k=0$ term of that sum simplifies to be the gamma distribution function part of the GPC convolution. Splitting off that term and adjusting the lower summation index from $k=0$ to $k=1$ yields,
 
\begin{equation}\label{eq:few2}
\frac{b^a }{\Gamma (a)} e^{-b\, t} t^{a-1} -\frac{ \alpha\, b^a }{\Gamma (a)}t^{a-1}\sum _{k=1}^{\infty } \frac{\left(\frac{\beta }{t}\right)^k}{k! (k-\alpha )}(1-a)_k \, _1F_1(a,a-k;-b\, t)\;.
\end{equation}

Next, the quickly convergent sum term, Eq.~\eqref{eq:few2}, is added to the gamma distribution plus asymptotic formula Eq.~\eqref{eq:asy} to create a series accelerated algorithm rewrite of Eq.~\eqref{eq:GPC} for long $t$-values, 

\[\begin{split}
\textnormal{GPC}\arraycolsep=1.2pt\def\arraystretch{.7}
\left(\begin{array}{cc}
a&b\\
\alpha&\beta
\end{array} \Big|\,t\right)=&-\theta (t-\beta) \frac{ \alpha\, b^a }{\Gamma (a)}t^{a-1}\sum _{k=1}^{\infty }\left(\mfrac{\beta }{t}\right)^k \frac{(1-a)_k}{k! (k-\alpha )} \, _1F_1(a,a-k;-b t)\\
&+\theta (t-\beta )\left[ \frac{b^a }{\Gamma (a)} e^{-b\, t} t^{a-1}-\pi \csc (\pi \alpha )\frac{b^a \beta ^{\alpha } }{\Gamma (\alpha )} t^{a-\alpha -1} \, _1\tilde{F}_1(a,a-\alpha ;-b t)\right]
\end{split}\;\;\;.\]

\noindent This is identically Eq.~\eqref{eq:accelgpc}, which completes the proof of the long-$t$ theorem. \end{proof}

The second line of the above equation is an asymptote of the GPC function. The above equation's first line when written as a list of terms to be summed has all negative elements when $k>\alpha$, which was the case for metformin \cite{Wesolowski_2020}. If $k<\alpha$ for the first few $k$, then the simplified summation terms are initially positive until $k>\alpha$, but in any case the magnitude of those terms is strictly monotonically decreasing such that increasing precision to sum those terms is unnecessary. The confluent hypergeometric functions in those terms and their effects on convergence are presented in detail in the \nameref{long} Appendix Subsection, which shows that the absolute value of the ratio of the $(k+1)$th to $k$th terms is approximately $\frac{\beta}{k\,t}$, where the $k$ in the denominator insures that the absolute values of the simplified terms of the summand for the above formula are monotonically decreasing, and that each $(k+1)^{\text{st}}$ term is many times closer to zero than the $k^{\text{th}}$ term, such that it is unnecessary to test for convergence using the sum to infinity of all the remainder terms, i.e., in practice it is sufficient to test the absolute value of the last term and to stop the summation when that magnitude is less than the desired precision (e.g., $<10^{-65}$).

\subsection*{Other long-\textit{t} functions; the integrals and derivative}
\textbf{GPC type I long-\textit{t} integral:} The derivation of a similarly accelerated series for $t\geq4\,\beta$ of the CDF of GPC, i.e., its $0\text{ to } t$ integral, GPC$_F$, follows from its primary definition, Eq.~\eqref{eq:CDF}, using the same procedure as Eqs.~\eqref{eq:GPC2} to \eqref{eq:accelgpc}, leading to, 

\begin{equation}\label{eq:fastgpc}
\begin{split}
\textnormal{GPC}_F\arraycolsep=1.2pt\def\arraystretch{.7}
\left(\begin{array}{cc}
a&b\\
\alpha&\beta
\end{array} \Big|\,t\right)=&-\theta (t-\beta) \frac{ \alpha\, b^a }{\Gamma (1+a)}t^{a}\sum _{k=1}^{\infty }\left(\mfrac{\beta }{t}\right)^k \frac{(-a)_k}{k! (k-\alpha )} \, _1F_1(a,a-k+1;-b\, t)\\
&+\theta (t-\beta ) \left[1-Q(a,b\,t)-\pi \csc (\pi\, \alpha )\frac{ b^a \beta ^{\alpha } }{\Gamma (\alpha )} t^{a-\alpha}\, _1\tilde{F}_1(a,a-\alpha +1;-b\, t)\right]
\end{split}\;\;\;,
\end{equation}

\noindent where $Q(a,b\,t)=\frac{\Gamma(a,b\,t)}{\Gamma(a)}$ is the regularised upper incomplete gamma function, and is the complementary cumulative density function (CCDF$=1-$CDF) of the gamma distribution.\footnote{CCDF is sometimes loosely referred to as a survival function, $S(t)$.} Note that GPC$_F$ is a CDF, such that the upper limit of Eq.~\eqref{eq:fastgpc} as $t$ increases is 1 or 100\% of initial dose eliminated from the body. 
\newline
\newline
\textbf{GPC type I long-\textit{t} double integral:} Similarly, the super cumulative distribution, i.e., the integral from $\tau=0\text{ to }t$ of the CDF is,

\begin{equation}
\begin{split}
\textnormal{GPC}_\mathcal{F}&\arraycolsep=1.2pt\def\arraystretch{.7}
\left(\begin{array}{cc}
a&b\\
\alpha&\beta
\end{array} \Big|\,t\right)=-\theta (t-\beta)\frac{\alpha\, b^a t^{a+1} }{\Gamma (a+1)}\sum _{k=2}^{\infty}\left(\frac{\beta }{t}\right)^k \frac{(-a)_k }{k! (k-\alpha )(a-k+1)} \, _1F_1(a,a-k+2;-b \,t)\\
&+\theta (t-\beta ) \bigg\{\frac{t\, e^{-b \,t} (b \,t)^a}{\Gamma (a+1)}-\frac{\alpha\, \beta }{\alpha -1}[1-Q(a,b \,t)]+\left(t-\frac{a}{b}\right) [1-Q(a+1,b \,t)]\\
&\hspace{5em}-\pi\, \csc (\pi\, \alpha )\frac{b^a \beta^\alpha }{\Gamma (\alpha )}\,t^{a-\alpha+1 }\, _1\tilde{F}_1(a,a-\alpha +2;-b \,t)\bigg\}
\end{split}\;\;\;.
\end{equation}

\noindent Note that the sum term is now indexed from $k=2$, for which each simplified summation element has a negative value when $k>\alpha$, and a multiplied out positive first term when $\alpha<2$. 
\newline
\newline
\textbf{GPC type I long-\textit{t} derivative:} The GPC derivative's algorithm for $t>4\beta$, i.e., long-$t$, is
\begin{equation}
\begin{split}
\textnormal{GPC}\,'&\arraycolsep=1.2pt\def\arraystretch{.7}
\left(\begin{array}{cc}
a&b\\
\alpha&\beta
\end{array} \Big|\,t\right)=\theta (t-\beta)b^a t^{a-2}\Bigg\{\frac{ a-b\,t-1}{\Gamma (a)}e^{-b\,t}\\
&+ \frac{\pi \alpha \csc (\pi\, \alpha ) \left(\frac{\beta }{t}\right)^{\alpha } }{\Gamma (\alpha +1)}\left[(\alpha +1) \, _1\tilde{F}_1(a;a-\alpha ;-b \,t)-a \, _1\tilde{F}_1(a+1;a-\alpha ;-b\, t)\right]\\
&+\frac{\alpha }{\Gamma (a)} \sum _{k=1}^{\infty } \frac{(1-a)_k \left(\frac{\beta }{t}\right)^k }{k! (k-\alpha )}\left[b\, t \, _1F_1(a;a-k;-b\, t)-(a-k-1) \, _1F_1(a-1;a-k-1;-b\, t)\right]
\Bigg\}
\end{split}\;.
\end{equation}

\subsection*{The combined short- and long-\textit{t} algorithm for GPC series acceleration}
There are now two algorithms, an algorithm that converges quickly only for short $t$-values, and another that converges quickly only when $t$-values are long. This section describes how the algorithms are combined to produce a new accelerated algorithm for any value of $t$. A full set of functions for the derivative and integrals of the GPC algorithm follows the same pattern as the \nameref{alg} Appendix Subsection. The two algorithms are combined by choosing $t=4\,\beta$ as the floor (least) value for use of the long-$t$ algorithm, makes the next term at worst approximately 1/4 that of the current term. Given a next term fraction of $\frac{\beta}{k\,t}$ times the current term, the $t=4\,\beta$ floor value is not critical, the trick is to avoid second to first term ratios that initially approach 1 as $t\to\beta$, for which the short-$t$ algorithm has fewer terms and converges faster. See the \nameref{choose} Appendix Subsection for further information.

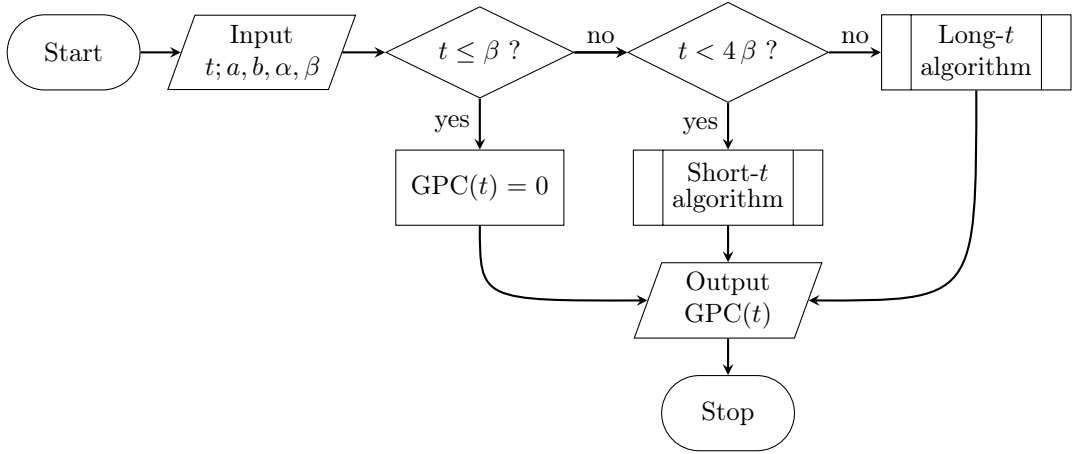
\begin{figure}[ht] \centering
\begin{tikzpicture}[node distance=1.5cm]
\node (start) [beginend, xshift=2cm] {Start};
\node (in1) [io, right of=start, xshift=1cm] {Input $t;a,b,\alpha,\beta$};
\node (dec1) [decision, right of=in1, xshift=1.4cm] {$t\leq\beta$ ?};
\node (dec2) [decision, right of=dec1, xshift=1.8cm] {$t<4\,\beta$ ?};
\node (pro2b) [process, below of=dec1, yshift=-.3cm] {GPC$(t)=0$};
\node (pre1) [draw,rectangle split, rectangle split horizontal,rectangle split parts=3,minimum height=1cm, right of =dec2, xshift=1.8cm] {\nodepart{two}\shortstack{Long-$t$\\algorithm}};
\node (pre2) [draw,rectangle split, rectangle split horizontal,rectangle split parts=3,minimum height=1cm, below of =dec2,yshift=-.3cm] {\nodepart{two}\shortstack{Short-$t$\\algorithm}};
\node (out1) [io, below of=pre2] {Output GPC$(t)$};
\node (stop) [beginend, below of=out1] {Stop};
\draw [arrow] (start) -- (in1);
\draw [arrow] (in1) -- (dec1);
\draw [arrow] (dec1) -- (dec2);
\draw [arrow] (dec1) -- node[anchor=south] {no} (dec2);
\draw [arrow] (dec1) -- node[anchor=east] {yes} (pro2b);
\draw [arrow] (dec2) -- node[anchor=south] {no} (pre1);
\draw [arrow] (dec2) -- node[anchor=east] {yes} (pre2);
\draw [arrow] (pre2) -- (out1);
\draw [arrow] (pre1.south) .. controls +(down:2.8cm) .. (out1.east);
\draw [arrow] (pro2b.south) .. controls +(down:1.cm) .. (out1.west);
\draw [arrow] (out1) -- (stop);
\end{tikzpicture}
\caption{Standard flow chart for the \nameref{alg} Appendix Subsection. The predefined short- and long-$t$ routines are also described in the text.} \label{tikzchart}
\end{figure} 

The program uses so-called infinite magnitude numbers such that numbers like $\pm10^{\pm100000}$ can be used without overflow or underflow (code: \$MaxExtraPrecision = $\infty$). However, there is another concern; precision. Machine precision was 53 bits, or approximately 16 significant figures. When $10^{-100}$ and 1 are added, one has to have a precision of 100 significant figures to avoid truncation. For the short-$t$ algorithm the extended precision needed is precalculated using machine precision of large numbers, which are stored as simplified terms, and are searched to find the largest magnitude number (code: Ordering[storage,$-1] [[1]]-1$). It is then that number as a rounded base 10 exponent (code: Round[Log10[Abs[outmax]]]) plus 65 significant figures that is used as the required precision of the computation. The terms of the summand are then recalculated to that high precision, then summed, such that the result has approximately 65 significant figures remaining even though the calculation itself may have needed a thousand or more significant figures to yield that result. The same approach could be used to calculate $\sin(x)$ from its infinite series definition. As mentioned above, in practice that is not used, and instead the equivalent principal sine values of $x$ are computed. For the GPC$(t)$ computation, one can invert the range related extra precision problem by reordering the series to make it increasingly less demanding to calculate long-$t$ values by direct application of Eq.~\eqref{eq:accelgpc} and that is precisely what the long-$t$ GPC type I algorithms does. The value $t=4\beta$ is used to transition between shorter $t$-values for use by the short-$t$ algorithm, and longer $t$-values for use with the long-$t$ algorithm. As mentioned, that time of transition between long and short algorithms is not critical and is more formally presented in the \nameref{choose} Appendix Subsection. 

\section*{Results}

This Results section shows examples for GPC algorithm run times and diagnostics, of how it can and should be used including the use for extended same dose multidosing, and a subsection illustrating confidence interval (CI) and coefficient of variation (CV) diagnostic quality assurance. 

\subsection*{Algorithm run time analysis}
\begin{figure}[H]
 \centering
 \includegraphics[scale=.325]{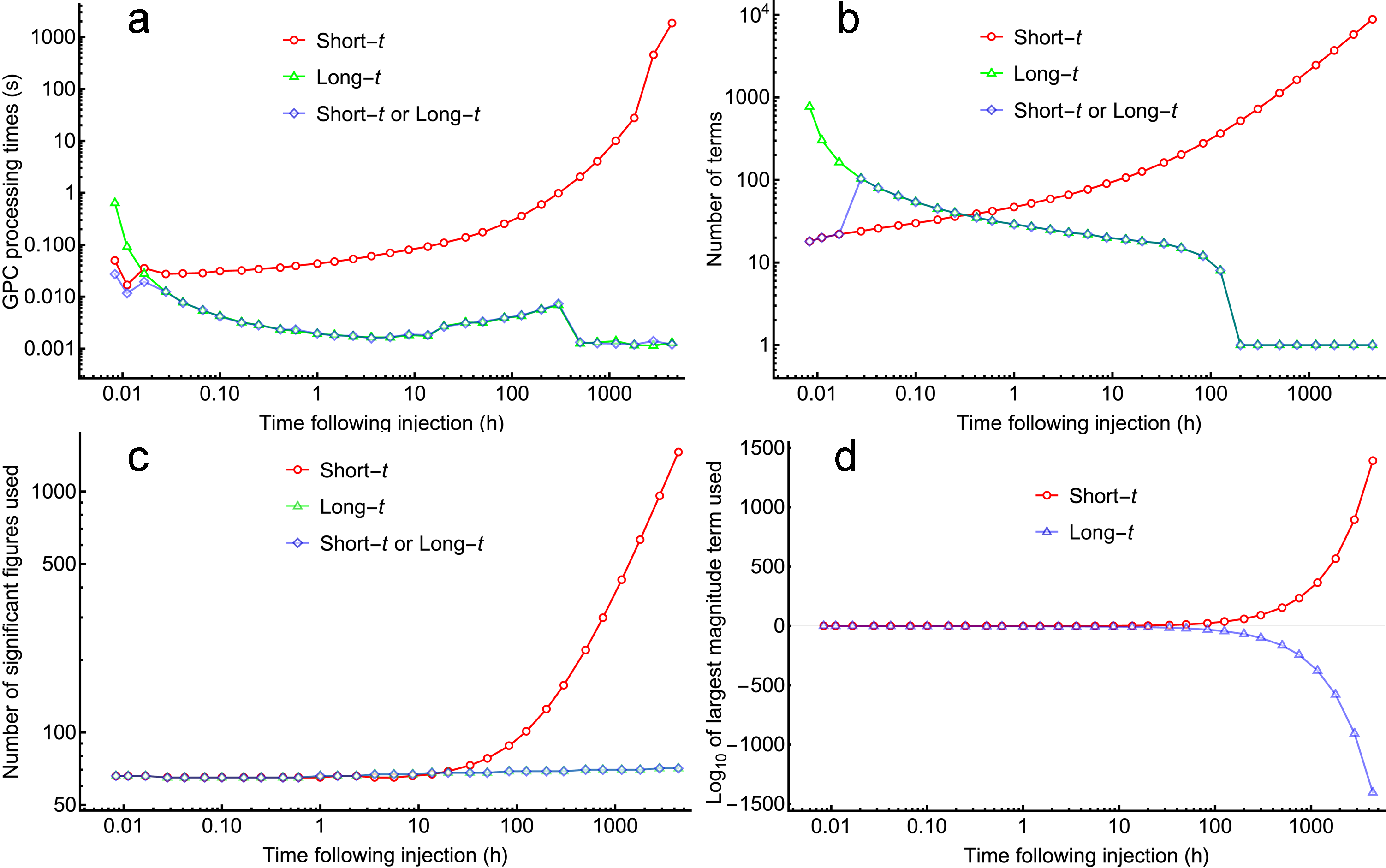} 
 \caption {Log-log plots comparing the performance of the short time (red connected open circles), the long time algorithm (green connected open triangles) and the either short or long time (blue connected open diamonds) algorithms' performance for the GPC model of metformin disposition in dog 1 of reference \cite{Wesolowski_2020}. Panel \textbf{a} shows that the long or short-$t$ algorithm is more than one million times faster than the short-$t$ algorithm when the time after injection is very long, e.g., for predicting a serum concentration at one half year (4396 h), and more than 20 times faster than the long-$t$ algorithm for $t=30$ s. Panel \textbf{b} shows the number of summed terms $(n)$ for each method. Note that for long $t$-values, the combined algorithm used the long-$t$ method and only calculates one term, whereas the short-$t$ algorithm would have used many terms. Panel \textbf{c} shows the precision needed to accommodate the largest of the terms for each algorithm. As $t$ increases, the short-$t$ function uses an increasing number of significant figures, but for the long-$t$ or combined algorithms that number increases only slightly to preserve accuracy for lower concentrations at longer times. Panel \textbf{d} shows the largest term magnitude as powers of 10 for the short- and long-$t$ algorithms. For long times, the short-$t$ algorithm alternating sign intermediate terms reach quite large magnitude, while for the long-$t$ algorithm the largest magnitude term collapses to vanishingly small values.}
 \label{Fig-2}
 \end{figure}

The GPC combined short- and long-$t$ algorithm was defined in terms of how to calculate it efficiently, as above. Implementation of the combined short and long time algorithm using Mathematica 12.3 without parallel processing on a 2.3 GHz 8-Core Intel i9 processor allows long $t$-value execution times of around 1.2 millisecond with typical 63 to 67 decimal place accuracy. (The full range of run times is approximately from 42 to 1.2 milliseconds for $t$-values ranging 30 s to 1/2 year.) This contrasts with the short-$t$ implementation of GPC Eq.~\eqref{eq:GPC2}, which, as $t$ increases, needs more terms and higher precision to maintain a given precision of the final result, with a processing time that progressively becomes intractably long. Figure \ref{Fig-2} shows the relative performance of the algorithms in 

\begin{figure}[H]
 \centering
 \includegraphics[scale=.32]{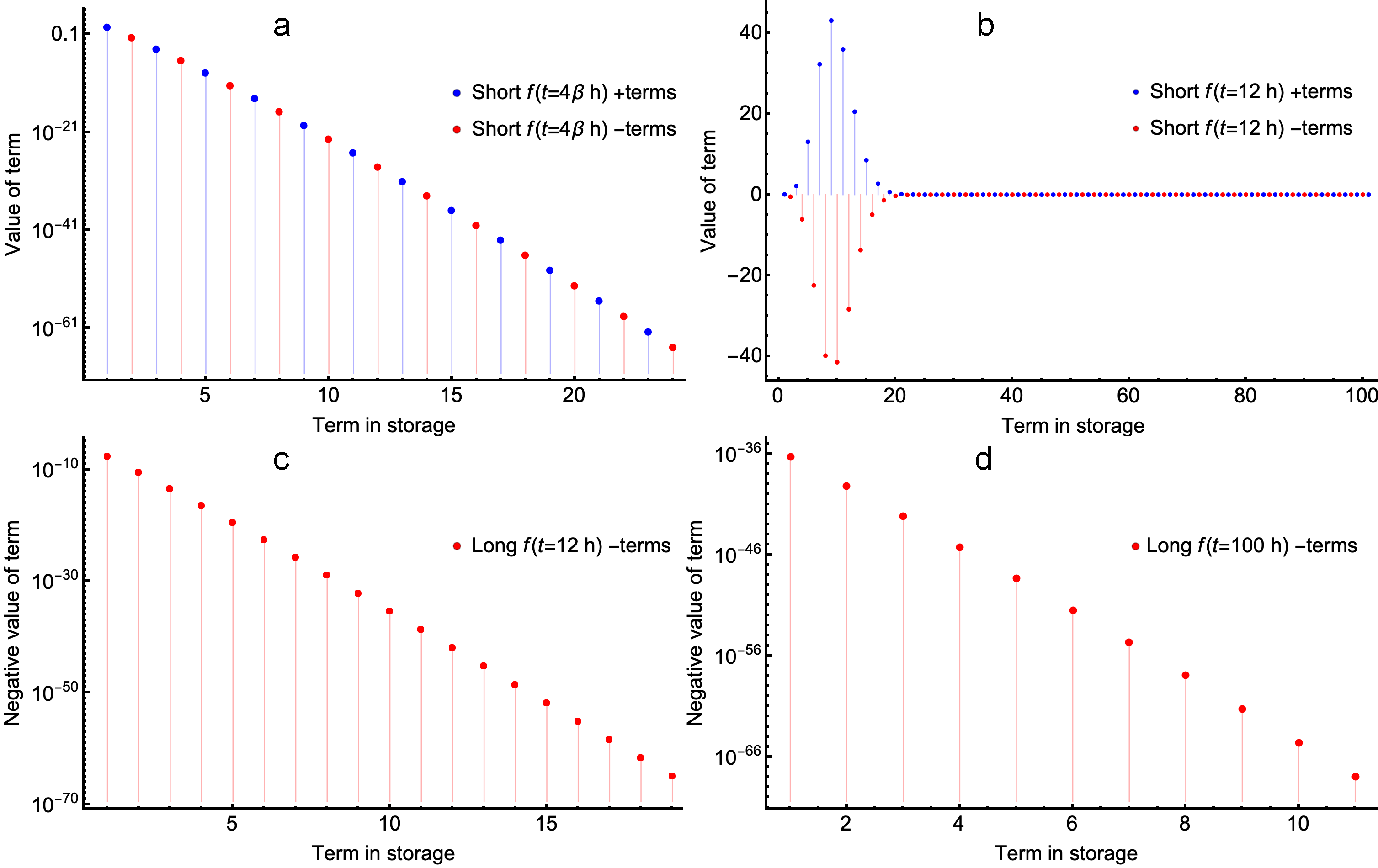} 
 \caption {Individual terms of the sums for the short-$t$ and long-$t$ algorithms for the GPC model of metformin disposition in dog 1 of reference \cite{Wesolowski_2020}. These terms are stored in a list variable called, somewhat unimaginatively, \textit{storage}. Note that $f(t=4\beta\, \text{h})$, is an explicit replacement rule meaning that \textit{an $f(t)$ is evaluated $f(4\beta)$, where $t=4\beta$ h}. Panels \textbf{a} and \textbf{b} show the values of the short-$t$ algorithm's summand terms. Panel \textbf{a} shows the values for the short-$t$ algorithm at the upper limit of $t$ of its usage in the combined, new, algorithm at $t=4\beta$ (100 s in this case). The blue dots are positive values, and the red dots are negative values of the summand of Eq.~\eqref{eq:GPC}. Panel \textbf{b} shows what happens when the short-$t$ algorithm is used at 12 h. That is the oscillatory terms would have intermediate values that grow in magnitude before they converge. While this poses little problem at 12 h, this time is not used in the new, combined algorithm and in the extreme the oscillatory intermediate terms of the summand grow to very large magnitude as $t$-values increase. Thus, before the summand is calculated for long $t$-values when using the short-$t$ algorithm, preliminary calculation of the required extended precision is necessary, which prologues execution time markedly. Panel \textbf{c} shows the region of values for the long-$t$ algorithm at 12 h. Note that there are fewer terms than for the short-$t$ algorithm at that same time (panel \textbf{b}), that all terms are negative (red), and that they decrease by one or more orders of magnitude between successive terms. Panel \textbf{d} shows that by 100 h for the long-$t$ algorithm, there are fewer terms than at 12 h, and that their magnitude is very small even for the largest magnitude, first, term.}
 \label{Fig-3}
 \end{figure}

\noindent these respects using the GPC parameters from fitting metformin data for dog 1 \cite{Wesolowski_2020}. This dog showed the median regression error of 8.7\% of the seven studied. Despite having the fastest elimination at 72 h, the concentration level for that dog was predicted to be $2 \times 10^{-7}$ of peak at one year, a small number but much larger than could be produced assuming a terminal exponential tail. For the short-$t$ algorithm the run-time to calculate concentration at one-half year following injection was 1809 s, versus 1.2 milliseconds for the new algorithm. This difference is because the short-$t$ algorithm used at long times had 8883 terms to sum, and the call to \textsf{gpcshort} was used twice; once at machine precision to find the maximum absolute value term $(1.0796*10^{1392})$ of all of the summand terms in order to calculate that 1457 place precision was required to obtain 65 place precision in the output, and once again to do the 1457 place summation arithmetic. For the combined (new) algorithm this is not needed as for short times the short-$t$ algorithm does not have large oscillating terms, and the long-$t$ algorithm has monotonically decreasing term magnitude both for each sequentially summed term, and as $t$ increases, for each first term magnitude. For example, the first (and only) term of the long-$t$ algorithm's summand at one-half year was negligible $(-1.851*10^{-1403})$. These effects are illustrated in Figure \ref{Fig-3}.

\noindent For our test case example, the two algorithms, short-$t$ and long-$t$, agreed to within 63 to 67 decimal places. In practice, the short-$t$ algorithm is used for short times and the long-$t$ algorithm is used for long times. It makes little difference what cut-point time between short- and long-$t$ algorithms is used, and the time $4\beta$, albeit around 100-120 s, was chosen as a division point between algorithms short enough to ensure that extra precision padding for the short-$t$ algorithm would be unnecessary. 
 
\subsection*{Regression processing elapsed times and extended multidosing} 
\begin{figure}[H]
 \centering
 \includegraphics[scale=.35]{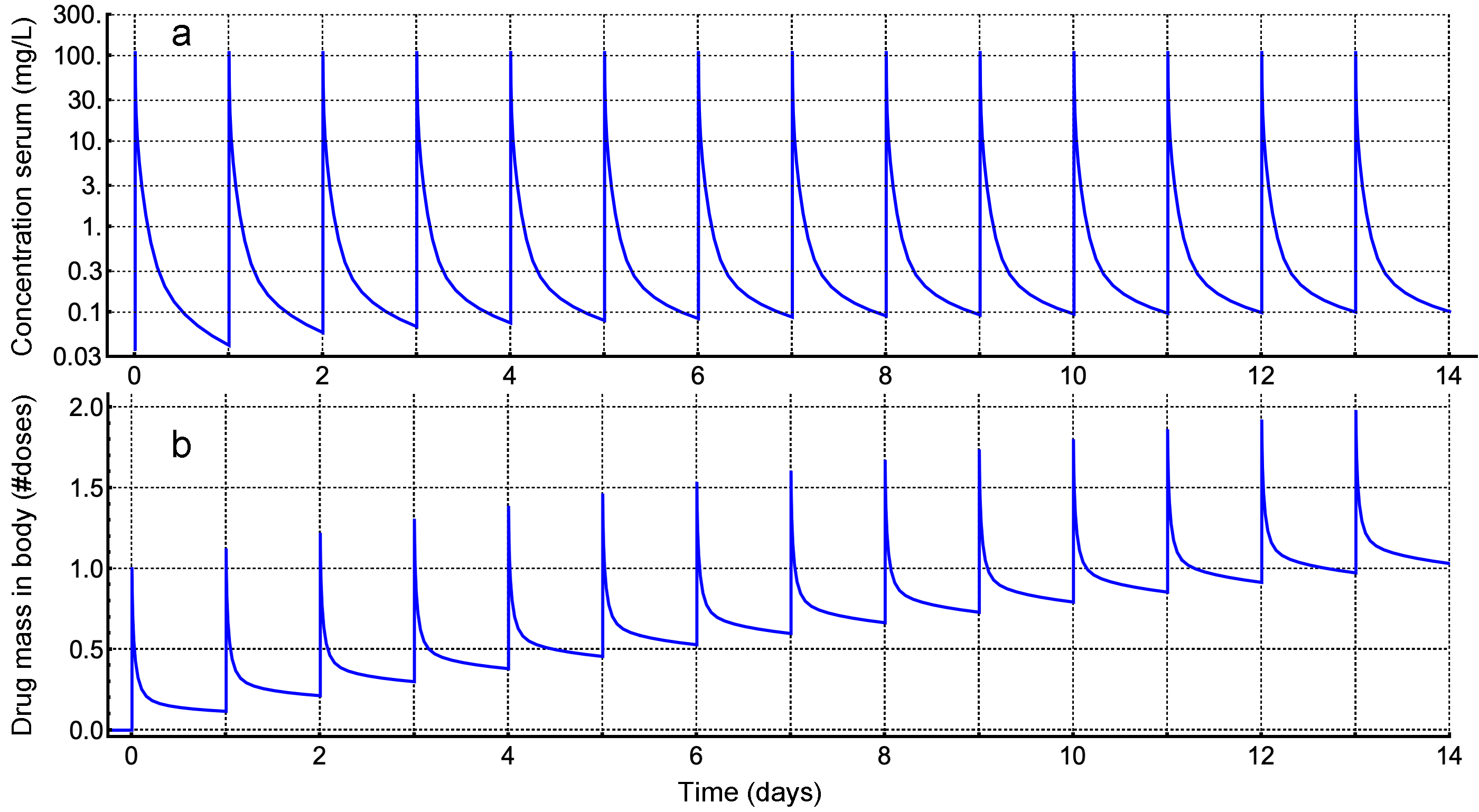} 
 \caption {Multidosing of the GPC function for dog 1 with one 18.248 mg/kg body weight dose every 24 h. Panel \textbf{a} shows the predicted concentration curve. Note that the peak concentrations did not increase much following 14 doses (0.089\%), but that the trough values increased rather more substantially (2.48 times) For panel \textbf{b} showing number of doses retained in the body, one first sees a 1 dose peak increasing to a 1.97 dose peak for the 14$^\text{th}$ dose, whereas the trough initially at 0.117 doses, increased 8.85 fold to finish at a 1.03 dose trough.}
 \label{multi_14-days}
 \end{figure}

For evaluating the 72 h data for seven dogs, the new, combined short- and long-$t$ algorithm run time for curve fitting was approximately 1:15 to 3:00 (min:s) average values, the program prior version with hardware and software accelerations for the short-$t$ algorithm and without sufficiently extended precision (despite using at least 65 place arithmetic) had run times in the approximate range of 34 to 35 min, but with occasional errors in parameter values of $\pm2\times10^{-14}$. With proper precision extension the error dropped below $10^{-20}$ for all 5 parameters and 7 cases, but the run time increased to 50 min, using a partly accelerated short-$t$ algorithm (Eq.~\eqref{eq:proved}) and 8-core hardware acceleration. The current combined short- and long-$t$ algorithm does not use those additional accelerations. Forty model-based bootstrap cases generated for the first dog's data---see next Subsection---took 49:45 (min:s), or 1:15 per case. That is a lot faster than the 33:51 per case it took to generate 35 bootstrap models using the old software (19:44:55). Overall, the run time is very approximately 27 times faster than prior, but is variable depending on the problem being solved, the computer hardware used, and the other software running on the computer at that time. For example, Figure \ref{multi_14-days}a, with a current run time of 7.1 s, could not be computed at all using the earlier software version.

\noindent Notice that if we wish to glean information during metformin multidosing with plasma or serum sampling, the best time to do so is just prior to the next scheduled dosing as those concentrations change for each dose interval, whereas the peak concentration change over time is very small. However, because the tissue dosage\footnote{For a single dose, body drug mass is $M(t) = \text{Dose } [1-\text{GPC}_F(t)]$} accumulates, the amount of drug in the body (Figure \ref{multi_14-days}b) cannot be predicted from serum (or plasma) concentration alone. Note that approximately one entire dose has accumulated in tissue by 14 days despite most of the cumulative dosage having been eliminated over that time. That is, during the first dose interval, the mean drug mass remaining in the body was 0.175 doses, and during the 14$^\text{th}$ dose interval the mean drug mass remaining in the body was 1.118 doses, where 12.88 dose masses were eliminated. 

\subsection*{Which are better, confidence intervals or coefficients of variation?}With reference to Table \ref{CIQ}, confidence intervals (CI) of the mean were extracted from model-based bootstrap with 40 cases generated for the first dog's data. For calculating CI's of the mean, the Student's-$t$ method was used (Verified assumptions: Central Limit Theorem, $n>30$, light-tailed distributions). However, as a result of extensive testing the degrees of freedom were set at $n$ rather than the more typical $n-1$, as it was found that for smaller values of $n$, physically impossible results were obtained, whereas even for $n=2$, when $n$ was used, rather than $n-1$, the results were accurate. For $n=40$ it made very little difference whether $n-1$ or $n$ were used. Also shown are CI's of the model based  bootstrap (A.K.A., parametric bootstrap) results calculated directly from the $n=40$ data using the nonparametric quantile (A.K.A, percentile) method of Weibull \cite{gurland1971simple}.\footnote{This uses the Weibull method for extracting confidence intervals, which in Microsoft Excel (2007) would format for the lower tail as PERCENTILE.EXC(A1:A40,0.025) and from Mathematica 12.3 \cite{Mathematica} as Quantile[data, 0.025, $\{\{0, 1\}, \{0, 1\}\}$], https://mathworld.wolfram.com/Quantile.html} Note that the Pareto rate parameter, $\beta$ was not presented. Since many (38 of 40) of the results were are the constraint boundaries of 25 to 30 s, one already knows what the confidence interval largely is; the constraint values themselves. Another situation entirely exists for coefficients of variation (CV). Note in the table that when $n=5$ as during the prior study, that the values so obtained were too small. It is theoretically possible to use bootstrap (in our case that would be bootstrap of model-based bootstrap) to obtain confidence interval quantiles for the median CV, and although median values of CV's have shown sufficient robustness to construct confidence intervals for $n$ sufficiently large \cite{brody2002significance}, the correction for $n$-small is problematic as per Table \ref{CIQ} and the \nameref{Discussion} Section that follows. 

 \begin{table}[ht]
\centering
\captionsetup{justification=justified,margin=0cm}
 \caption{Example parameter results and quality assurance are shown for dog 1. The relative root mean square error of fit (rrms \%) and R$^2$ values are acceptable. Confidence intervals (CI) are shown for the parameter values and the mean parameter values.} 
 \label{CIQ}
\vspace{-.5em}
\setlength\tabcolsep{2.7pt}
\begin{tabularx}{.8\columnwidth}{cccccccc}
\hline

\vspace{-0.2em}&Primary&95\% CI&Mean&95\% CI&SD&CV\%&Prior\\
&result&bootstrap&&of mean&&&CV\%\\
\Xhline{2\arrayrulewidth}
$n$&1&40&40&40&40&40&5\\
rrms \%&8.71&4.72 to 9.88&7.52&7.16 to 7.89&1.15&15.3&---\\
R$^2$&0.99872&0.99773 to 0.99950&0.99872&0.99860 to 0.99884&0.00038&0.039&---\\
$a$&0.349&0.198 to 0.516&0.350&0.324 to 0.376&0.0802&22.9&8.09\\
$b$&0.732&0.558 to 0.889&0.735&0.708 to 0.761&0.0821&11.2&4.84\\
$\alpha$&0.264&0.235 to 0.322&0.268&0.261 to 0.275&0.021&7.85&5.13\\
\textit{AUC}&31.2&26.8 to 38.0&31.3&30.4 to 32.2&2.77&8.84&6.22\\
\textit{CL}&9.76&8.01 to 11.3&9.78&$\;\,$9.52 to 10.05&0.829&8.47&6.42\\
\hline
\end{tabularx}
\end{table}

\section*{Discussion}\label{Discussion}

Wise \cite{Wise1985} first proposed that power functions or gamma distributions should be used in pharmacokinetic modelling as superior alternatives to sums of exponential terms. This has been reinforced more recently, for example by Macheras \cite{Dokoumetzidis2010}. While convolution models and fractal consistent models have been shown to be superior models in some cases and find occasional use \cite{garrett1994bateman,Wesolowski2016GDC,wanasundara2016,Wesolowski2016PLoS,Wesolowski_2020} compartmental modelling software is widely available and is used by default. For example, compared to biexponential clearance evaluation of 412 human studies using a glomerular filtration rate agent, adaptively regularised gamma distribution (Tk-GV method \cite{wesolowski2010tikhonov,wesolowski2014method}) testing was able to reduce sampling from 8 to 4 h and from nine to four samples for a more precise and more accurate, yet more easily tolerated and simplified clearance test \cite{wanasundara2016}. Despite this, few institutions have implemented the Tk-GV method at present. In the case of metformin, a highly polar ionised base, the extensive, obligatory active transport of the drug into tissue produces a rate of expansion of the apparent volume of distribution having the same units as renal clearance, yielding the Pareto (power function) tail. This heavy tail, and Figure \ref{multi_14-days}, help to explain why metformin control of plasma glucose levels had delayed onset, e.g., following 4-weeks of oral dosing \cite{buse2016primary}, and provides hints concerning the lack of a direct correlation between drug effect and blood metformin concentrations \cite{stepensky2002pharmacokinetic}. Other basic drugs whose active transport dominates their disposition may show similar behaviour. The long tail in the disposition of amiodarone may be a reflection of its very high lipid solubility rather than, or in association with, active tissue uptake. Weiss \cite{Weiss1999} described its kinetics after a 10 min intravenous infusion with an \textit{s}-space Laplace transform convolution of a monoexponential cumulative distribution with an inverse Gaussian distribution and a Pareto type I density, which lacked a real or \textit{t}-space inverse transform such that the modelling information had to be extracted numerically. A real space $f(t)$ model convolution of time-limited infusion effects of a GPC type I distribution is simple to construct and would be the same as Weiss's model in the one essential aspect that matters; testing of the amiodarone power function tail hypothesis, for which a GPC derived model would have the advantage of being more transparently inspectable. Similarly, Claret \textit{et al.} \cite{Claret2001} used finite time difference power functions to investigate cyclosporin kinetics for which GPC and related model testing may be appropriate. 


We were able to use Nelder-Mead global search regression model-based bootstrap to provide more information and better information for parameter variability than would be available from a gradient matrix. Some readers would prefer to use the Levenberg-Marquardt algorithm convex gradient regression method, so that the gradients can be used to estimate case-wise coefficients of variation. The logarithm of sums of exponential terms is always convex. The GPC-metformin loss function is nonconvex, as shown by failure of the interior point method to improve on solutions as reported in the \nameref{data} Subsection. Constrained nonconvex gradient methods are comparatively rarely implemented; there appears to be no such implementation in Mathematica at present. 

Correction of standard deviation (SD) for small numbers ($n<30$) using bootstrap of model-based bootstrap and $\chi^2$ were used as mentioned elsewhere \cite{friedman2009elements}, and led to using $n$ rather than $n-1$ for Student's-$t$ degrees of freedom. Whereas variance is unbiased, when the square-root of variance is taken, the result, standard deviation becomes biased. Arising from $\chi^2$, a standard deviation from only two samples, is on average only 79.8\%, $\sqrt{\frac{2}{\pi }}$, of the population standard deviation \cite{gurland1971simple}.\footnote{\label{note1}Given only two samples, the population mean is not located midway between them, however, the midpoint (mean) is used to estimate the population mean in the standard deviation formula. The correction formula multiplier for an unbiased estimator ($\hat{\sigma}$) of population standard deviation ($\sigma$) from sample standard deviation ($s$) is $\hat{\sigma}=c_n s$, where $c_n=\sqrt{\frac{n-1}{2}} \Gamma \left(\frac{n-1}{2}\right)\Gamma \left(\frac{n}{2}\right)^{-1}$ \cite{gurland1971simple}.} Gradient methods lack \textit{pre hoc} testing of the implicit assumption of residual normality and do not \textit{post hoc} provide any parameter distribution information. From the trace of the gradient matrix, one obtains a standard deviation with degrees of freedom that are $n-p-1$ ($n$-samples, $p$ parameters) \cite{friedman2009elements}. For standard deviations in the case where $n-p-1$ is small, the corrections for standard deviation are large. Overall, the ratio between gradient based error propagation results and that from bootstrap is not unusually a factor of two larger or smaller \cite{green1987standard}. Moreover, average fit errors using any loss function >10\%, for assay methods with errors <10\% may suggest that the algorithm/data combination is suspect \cite{burger2010limited,Wesolowski2016GDC,zhang2016bootstrapping,Wesolowski_2020}, and for the metformin dog data that is the case for two- and three-compartment models, but not for the GPC model, which latter model was the only one to fit the data better than 10\% (average 8.6\% rrms with assay error of 5.2\% rrms), as well as being the only model to exhibit normality and homoscedasticity of residuals \cite{Wesolowski_2020}. When the fit error is >10\%, one should, at a minimum, test residuals for homoscedasticity and normality, and if these are not present, a better fit model should be sought for its own sake, and bootstraping becomes problematic \cite{zhang2016bootstrapping}. 

The use of coefficients of variation is sometimes problematic. Suppose that we have lots of data, but because $\text{CV}=\text{SD}/$mean, if by chance in a particular case especially if we have small $n$, some of the multiply generated mean values may approach zero, which injects some erratically high CV-values into a distribution of values. It is for that reason, numerically instability, that the more data one has, the worse the mean CV-value can be, with the solution being to first calculate many CV values, and then take their median value \cite{brody2002significance}. Even though the mean value may be not useful, the median may be, and confidence intervals for CV could be established using bootstrap quantiles, but not by using the gradient matrix approach because correction for $n$-small is problematic. That is, for mean values that can be rewritten as being proportional and having an established maximum range, e,g., Likert scale minus 1 variables, correcting CV underestimation for small values of $n$ is possible. However, if, as is the case here, there is no theoretical maximum CV, one needs to invent a correction based upon the observed confidence intervals of the mean \cite{smithson1982relative}, such that CI-values are unavoidable for determining the meaning of the preponderance of CV results. Finally, comparison for significant differences between parameters for one subject versus another are easy to construct using CI, but more difficult to obtain for CV. Thus, CV-values cited without explicit quality assurance should be regarded as qualitative results.

\subsection*{Limitations}
 
 A major deficiency of the first article that applied and compared the gamma-Pareto type I convolution (GPC) model to other models \cite{Wesolowski_2020} was the lack of an algorithm that could be used for independent investigation and/or for application to other drugs. The accelerated algorithm presented herein is the first publication of code for a gamma-Pareto type I convolution (GPC). As such, the algorithm was kept in a simple form without using all possible acceleration tools or stopping conditions. While it could be optimised for even shorter run-times using vector addition of subsets of the summands, by using Eq.~\eqref{eq:proved} to reduce summand magnitudes for the short-$t$ algorithm and/or combining partial sums of the summands for the short- or long-$t$ algorithms, by eliminating diagnostic parameters such as run-time calculations, by compiling it, and by multiple other means. However, that would be at the expense of clarity and/or simplicity of presentation. It is complicated to compute the values of functions like the $\sin(x)$ efficiently. For example, an answer with precalculated exact precision can be quickly generated for $\sin{x}$ using the CORDIC procedure, which is optimised for binary register operations at the machine level \cite{Volder1959}. At a much higher and slower level, compared to the GPC$(t)$ short-$t$ algorithm, the $\sin(t)$ function's series expansion has even larger magnitude terms for long $t$-values. In its current form, the combined short- and long-$t$ GPC algorithm is so much faster than the previously published run times using the seven dogs 72 h data and more generally valid that it is now a practical algorithm. The current implementation is no longer limited as to how long $t$ is, and the propagated error of up to $2\times10^{-14}$ for parameter values obtained from regression of 72 h data has been reduced to $<10^{-20}$. That error demonstrates the major extent to which errors from 65 decimal place precision can propagate during processing of tens of thousands of calculations, especially during regression, which typically, by default, halves the number of significant figures---see the \nameref{data} Subsection. This does not affect any of the parameter values listed in Table \ref{params}, but the ability to quickly calculate a larger number of model-based bootstrap results would improve the parameter CI estimates. Another consideration is how to obtain exactly $n$ significant figures precision when $n$ are requested. Currently, for 65 significant figures requested, a result precise to several significant figures greater or lesser than 65 is returned and the algorithm is written only for 65 significant figure precision. Generalising this to request to obtain an arbitrary specific precision for a GPC functional height awaits the next several generations of algorithmic refinement.
 
 \section*{Conclusions}
The new GPC type I algorithm consists of two parts, one for very early times, and another for late times. At times less than approximately $4\,\beta$, i.e., 100-120 s for the metformin data, the short-$t$ algorithm is actually faster than the long-$t$ algorithm. For early data, the short-$t$ algorithm has alternating sign terms of monotonically decreasing magnitude. However, when used at long times, the short-$t$ GPC algorithm required precalculation of the precision needed for later summation, which represents an improvement over the algorithm previously used \cite{Wesolowski_2020}. In the newly-proposed, combined short and long-$t$ algorithm this precalculation is unnecessary because of the long-$t$ algorithm usage for all but the shortest $t$-values, resulting in markedly accelerated convergence, and the new ability to predict concentration at any time, no matter how long.

\section*{Acknowledgements} 
The authors thank Kunal Khadke at Wolfram Research for assistance with precision and Mathematica block structures, and William J. Jusko of the University of Buffalo for his generous advice concerning the intellectual content.

\section*{Appendix}
This section provides information concerning convergence of the short- and long-$t$ algorithms, when they should be used, and how to encode them in the Mathematica \cite{Mathematica} language.
\subsection*{Short-\textit{t} GPC convergence}\label{short}
The short-$t$ algorithm is an alternating series sum. For alternating series one can distinguish two types of convergence. Conditional convergence in which the value of the infinite sum depends on the order in which summation is performed as shown by the Riemann rearrangement theorem, and absolute convergence for which any order, or permutation, of summation process yields the same, unique, sum. Convergence is defined as conditional when an alternating series converges but its absolute value does not \cite{agana2015classical}. For example, the alternating harmonic series $\sum_{n=1}^\infty \frac{(-1)^{n+1}}{n}$ has an absolute value ratio of next term to current term of $\frac{n}{n+1}$, whose limit as $n\to\infty$ is 1. That means that as $n$ increases, the next term approaches the same size as the $n$th term, such that the absolute sum of terms is not bounded above, and the order of addition of the original series determines what the total sum is, making changes in order of summation yield different, i.e., ambiguous, results. If a limiting term ratio is less than 1, for example $\frac{1}{2}$, the series is absolutely convergent, e.g., the limiting infinite sum ratio of $\frac{1}{2}$ for some eventual term is, in binary arithmetic, $0.111111\dots_2\to1_2=1$. It is fair to call series, whose limiting absolute value term ratio is 0, eventually-rapidly convergent. In the case of the short-$t$ algorithm the infinite sum of its absolute values is, for sufficiently large values of $t$, a very large number. However, for any real valued time, $t$, no matter how large, there is a real number, $M$, that is greater than the magnitude of the infinite sum of absolute values of Eq.~\eqref{eq:GPC2}. That $M$ can be fantastically large, but never infinite, makes it difficult to use Eq.~\eqref{eq:GPC2} without precision that is explicitly extended for the purpose of accurately forming the infinite sum for certain values of the parameters and long times, but it does not make convergence conditional. 

\newtheorem{2em}{Lemma}
\begin{2em}\label{Lemma-1}In the case of the short-t GPC algorithm, convergence is absolute and eventually rapid, such that the Riemann rearrangement theorem prohibition for resequencing infinite sums does not apply.\end{2em}

\begin{proof}To show absolute convergence of the short-$t$ GPC type I algorithm, we construct the absolute value of the ratio of the $(n+1)$th to $n$th term. First, we take the infinite series definition of the incomplete beta function, \footnote{http://functions.wolfram.com/06.19.06.0002.01}\[B_z(A,B)=z^A \sum _{j=0}^{\infty } \frac{z^j (1-B)_j}{j! (A+j)}=z^A\left[\frac{1}{A}-\frac{B - 1}{A + 1} z +\frac{(B-1) (B-2) }{2! (A+2)}z^2-\dots\right];\;\;| z| <1\land \neg (-a\in \mathbb{Z}\land -a\geq 0)\;.\]Although this is an alternating sign series with a restricted range of convergence, we term by term, \textit{without permutation}, substitute into it the incomplete beta function parameters of Eq.~\eqref{eq:GPC2}'s $n$th and $(n+1)$th terms; $B_{1-\frac{\beta}{t}}(a+n,-\alpha),$ and $B_{1-\frac{\beta}{t}}(a+n+1,-\alpha),$ and substitute that into the absolute value of the $(n+1)$th to $n$th term ratio of the summand of Eq.~\eqref{eq:GPC2}, and simplify to yield,\[\frac{b\, (t-\beta )}{n+1}\;\frac{\mfrac{1}{a+n+1}+\mfrac{(\alpha +1) }{a+n+2}\left(1-\mfrac{\beta }{t}\right)+\mfrac{(\alpha +1) (\alpha +2) }{2! (a+n+3)}\left(1-\mfrac{\beta }{t}\right)^2+\dots}{\mfrac{1}{a+n}+\mfrac{(\alpha +1) }{a+n+1}\left(1-\mfrac{\beta }{t}\right)+\mfrac{(\alpha +1) (\alpha +2)}{2! (a+n+2)} \left(1-\mfrac{\beta }{t}\right)^2+\dots}<\frac{b\, (t-\beta )}{n+1}\;\;.\]As $\infty>t>\beta$, neither the infinite series numerator or denominator is alternating, thus their ratio is absolutely convergent as $\mfrac{b\, (t-\beta )}{n+1}$ is an asymptote of, and upper bound for, the ratio of consecutive absolute value terms as $n\to\infty$. While $b\, (t-\beta )>n+1$, if that occurs, for example for long $t$-values, we would expect the magnitude of the terms of the summands to increase for $n$ small enough, but as $n$ increases $b\, (t-\beta )\ll n+1$ eventually, and the $(n+1)$th relative term magnitude can be made as asymptotically close to zero as desired, and convergent by the ratio test \cite{laugwitz1994riemann}. \end{proof}Thus, the magnitude of alternating terms is eventually monotonically decreasing such that the absolute error of summation from truncating at an $n$th term for $n$ sufficiently large is less than the magnitude of the $(n+1)$th term by the alternating series remainder theorem. Moreover, the first term of the summand is some definite positive real number proportional to 1. Setting the first term to be 1, we conclude that the sum of the absolute value of the summands of Eq.~\eqref{eq:GPC2} is proportional to a number bounded above by \[M\propto\sum _{k=0}^{\infty } \frac{[b (t-\beta )]^k}{k!}=e^{b (t-\beta )}\;\;,\]such that the sum of absolute values of summands of Eq.~\eqref{eq:GPC2} is bounded above by some positive constant value times an exponential function of $t$, and Eq.~\eqref{eq:GPC2} is absolutely convergent.

\subsection*{Long-\textit{t} GPC algorithm convergence rapidity}\label{long}
This subsection examines the rapidity of convergence of the long-$t$ GPC algorithm. In Lemma \ref{Lemma-1} directly above, it was shown that the short-$t$ algorithm is absolutely convergent. Therefore, its infinite series rewrite as the long-$t$ Theorem \ref{longT}, Eq.~\eqref{eq:accelgpc}, is also convergent but how many summation terms are needed for convergence and which parameters determine this convergence can be clarified using the substituted definition of the confluent hypergeometric series\footnote{http://functions.wolfram.com/07.20.06.0002.01} as follows.
\[\, _1F_1(a;a-k;-b\, t)=\sum _{j=0}^{\infty } \frac{(a)_j (-b\, t)^j}{j! (a-k)_j}=1-\frac{a\, b\, t}{a-k}+\frac{a (a+1) (b\, t)^2}{2! (a-k) (a-k+1)}-\dots\;\;.
\]
Note that in the limit as $k\to\infty$ the above equation is asymptotically ($\sim$) 1. Next, the ratio of the $(k+1)$th to $k$th term is asymptotic to $\mfrac{\beta }{k\, t}$ for $k$ sufficiently large,
\[\frac{\beta}{t}\; \frac{\alpha-k }{(k+1)( \alpha-k-1)}\;\frac{1-\mfrac{a }{a-k-1}b\, t+\mfrac{a (a+1) }{2 (a-k-1) (a-k)}(b\,t)^2-\cdots}{1-\mfrac{a}{a-k}b\, t+\mfrac{a (a+1) }{2 (a-k) (a-k+1)}(b\,t)^2-\cdots}\sim\frac{\beta }{k\, t}\;\;.\]


\noindent For that reason, for longer $t$-values, one can expect faster convergence of the long-$t$ algorithm with fewer terms summed. 
\subsection*{Choosing when to use the short-\textit{t} and long-\textit{t} algorithms}\label{choose}
\begin{figure}[H]
 \centering
 \includegraphics[scale=.55]{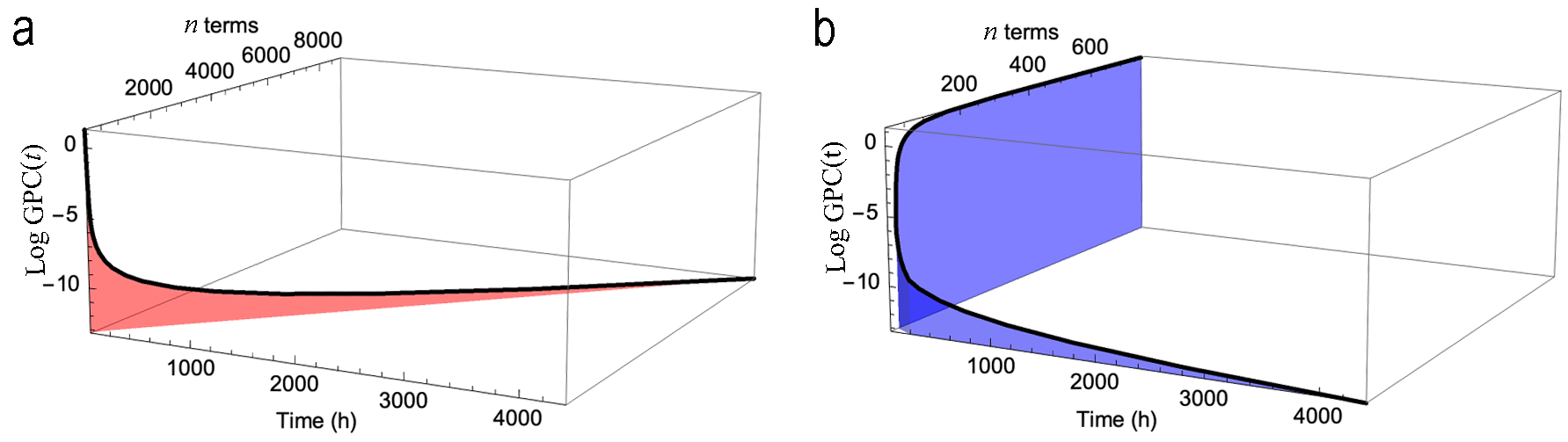} 
 \caption {The tradeoff between the time elapsed following bolus intravenous injection ($x$-axis), the number of terms to be summed for calculating the GPC function's value ($y$-axis), and the natural logarithm of the GPC function at that time ($z$-axis). Panel \textbf{a} shows short-$t$ algorithm performance and panel \textbf{b} shows the long-$t$ algorithm performance. Note that only for very early elapsed times does the short-$t$ algorithm have fewer terms than the long-$t$ algorithm.}
 \label{shortlong3d}
 \end{figure}

As above, the absolute value of the ratio of the next term to the current term for the short-$t$ algorithm is bounded above by $\mfrac{b\, (t-\beta )}{n+1}$. For the long-$t$ algorithm, the ratio of the $(k+1)$th to $k$th term approaches $\mfrac{\beta }{k\, t}$ for $k$ sufficiently large. Note that these are in the opposite direction, that is, while $t$-values increase, $\mfrac{b\, (t-\beta )}{n+1}$ increases and $\mfrac{\beta }{k\, t}$ decreases. It is not critical exactly at what $t$ value one elects to use the short- and long-$t$ algorithms, as the major cost in computational time and number of terms needed occurs at the extreme values of $t$, but in an opposite direction for each algorithm. Figure \ref{shortlong3d} shows the tradeoff for dog 1 of the metformin series between numbers of terms for summation, time following bolus injection, and the magnitude of the natural logarithm of $\text{GPC}(t)$, where GPC$(t)=\frac{C(t)}{\textit{AUC}}$. Selecting $t=4\beta$ as a cut point for switching between algorithms means that the short-$t$ algorithm absolute sum of terms is bounded above, from substitution into $e^{b (t-\beta )}$, by $e^{3b\,\beta}$ times the first term's value, not a large number, and the long-$t$ algorithm has an approximate maximum $(k+1)$th to $k$th term ratio of $\frac{1}{4k}$ for the shortest $t$-value used, which can still be made as small as desired for $k$ sufficiently large.

\subsection*{Mathematica source code of the GPC type I accelerated algorithm}\label{alg}
\textcolor{teal}{
\textsf{(*.....................The gpc$[t]$ function; a gamma-Pareto type I convolution fast calculation algorithm.............................\\
...........................................................Copyright Carl A. Wesolowski, 2021................................................................\\
To fit disposition data, minimise the loss function using AUC $\times$ gpc$[t]$ as the model, which returns AUC as a value. To use gpc$[t]$ enter the constrained > 0 coefficients a, b, $\alpha$, and $\beta$ to 65 decimal place accuracy. a and $\alpha$ are dimensionless, b is a rate and $\beta$ is a time. For example, the metformin dog 1 GPC type I parameters using b (h$^{-1}$) and $\beta$ (h) are*)}}\\
\textsf{a = N[Rationalize[0.34931003807815571524792421542558602868248355919027496611955665616], 65];\\
b = N[Rationalize[0.73182479199387479660419087183394451163091958778927254273673996698], 65];\\
$\alpha$ = N[Rationalize[0.26437129139517680335740710070693267536710608361890151476103695922], 65];\\
$\beta$ = N[Rationalize[0.0069444444444444444444444444444444444444444444444444444444444444444, 67];}\\
\textcolor{teal}{\textsf{(*Note the explicit input precision as 65 digits, which must be specified as such for the algorithm to function properly. During regression analysis the parameters values above would change to minimise a loss function, the ones above were used as realistic values for the algorithm execution timing trial of the Results Section.*)}}
{\setstretch{1.3}\\
\noindent \$\textsf{MaxExtraPrecision = }$\infty;$\\
\$\textsf{MinPrecision = 0;}\textcolor{teal}{\textsf{ (*This machine precision is adjusted in the Block commands below to $\$$MinPrecision = desirprec*)}}\\
\\
\textcolor{teal}{\textsf{(*.......................................Section for short times, $\beta<t<4\beta$, where $4\beta\approx $ 100 to 120 s.........................................*)}}\\
$\mathsf{gpcshort[\textcolor{blue}{\emph{t}}\_?NumericQ]:=Quiet\Big[Unevaluated\Big[k=0;storage=.;}$\textcolor{teal}{\textsf{ (*storage is an array for later summation*)}}\\
\indent $\mathsf{conscale= \frac{\alpha b^a \beta ^{\alpha }}{Gamma [a]}\textcolor{blue}{\emph{t}}^{\;a-\alpha -1};}$\textcolor{teal}{\textsf{ (*This is the constant multiplier of the sum and is calculated only once*)}}\\
\indent $\mathsf{target= \frac{10^{-65}}{ conscale};}$\\
\indent $\mathsf{storage=-conscale\;First@Last@Reap@While\left[Abs\left[Sow\big[ \frac{(-b \textcolor{blue}{\emph{t}})^k }{k!}Beta[1-\beta/\textcolor{blue}{\emph{t}}, a+k,-\alpha ]\big]\right]\right.}$\\
\indent \indent $\mathsf{ >target,k+\!+\Big];}$\\
\indent $\mathsf{nn= Ordering[storage,-1]\,[[1]]-1;}$ \textcolor{teal}{\textsf{ (*Finds nn, the index of the largest magnitude term in the storage array*)}}\\
\indent $\mathsf{outmax=conscale \frac{(-b \textcolor{blue}{\emph{t}})^{nn}}{nn!} Beta[1-\beta/\textcolor{blue}{\emph{t}},a+nn,-\alpha];}$\\
\indent $\mathsf{lastn=Length[storage];\;
xprec=Round[Log10[Abs[outmax]]];}$\\
\indent $\mathsf{If[xprec === Indeterminate, xprec = 0];}$ \textcolor{teal}{\textsf{ (*Exact times, e.g., $t=e^{1/1000}$ for $n=1$ may need this correction.*)}}\\
\indent $\mathsf{desirprec = 65 + Abs[xprec]\Big]\Big];}$\\
\\
 \textcolor{teal}{\textsf{(*...............................................$t\geq 4\beta$ section with asymptotes for long values of t...............................................*)}}\\
$\mathsf{gpcsetup[\textcolor{blue}{\emph{t}}\_?NumericQ]:=Quiet\Big[Unevaluated\Big[k=1;storage=.;}$\\
\indent $\mathsf{conscale = \frac{b^a \textcolor{blue}{\emph{t}}^{-1 + a}\alpha }{Gamma[a]};}$\textcolor{teal}{\textsf{ (*This constant multiplier is a different one than that above*)}}\\
\indent $\mathsf{target= \frac{10^{-65}}{conscale};}$\\
\indent $\mathsf{storage=-conscale\; First@Last@Reap@While\big[Abs\big[Sow[}$\\
\indent \indent $\mathsf{ \frac{(\beta/t)^k}{(k-\alpha) k!}Hypergeometric1F1[a,a-k,-b \,\textcolor{blue}{\emph{t}}\,]\; Pochhammer[1-a,k]]]>target,k+\!\!\,+];}$\\
\indent $\mathsf{nn = Ordering[Abs[storage], -1]\,[[1]];}$\\
\indent $\mathsf{outmax = -conscale \frac{(\beta/\textcolor{blue}{\emph{t}})^{nn}}{(nn - \alpha) nn!}
 Hypergeometric1F1[a, a - nn, -b\, \textcolor{blue}{\emph{t}}\,] \;Pochhammer[1 - a, nn];}$\\
\indent $\mathsf{lastn = Length[storage];}$\\
\indent $\mathsf{asympt[\textcolor{blue}{\emph{z}}\_?NumericQ] := \frac{b^a e^{-b \,\textcolor{blue}{\emph{z}}} \textcolor{blue}{\emph{z}}^{-1 + a}}{Gamma[a] }- \frac{\pi b^a \beta^\alpha \alpha Csc[\pi\,\alpha]}{Gamma[1 + \alpha]}\textcolor{blue}{\emph{z}}^{-1 + a - \alpha}Hypergeometric1F1Regularized[a, a - \alpha, -b \,\textcolor{blue}{\emph{z}}\,];}$\\
\indent $\mathsf{SetAttributes[asympt, Listable];
xprec = Round[Log10[Abs[outmax + asympt[\textcolor{blue}{\emph{t}}\,]\,]]];}$\\
\indent $\mathsf{desirprec = 65 + Abs[xprec]\Big]\Big];}$\\
\\
\noindent \textcolor{teal}{\textsf{(*..........................................GPC combined short and long time function call (use this one)...................................*)}}\\
$\mathsf{gpc[\textcolor{blue}{\emph{t}}\_?NumericQ] := Quiet\Big[If\Big[\textcolor{blue}{\emph{t}} \leq \beta, 0, If\big[\textcolor{blue}{\emph{t}} < 4 \beta, }$\\
\indent $\mathsf{Block\big[\{{\textcolor{blue}{\$MinPrecision}}= desirprec, {\textcolor{blue}{\$MaxExtraPrecision}} = \infty\},gpcshort[\textcolor{blue}{\emph{t}}\,]; \sum_{j=1}^{j=lastn}storage[[\,j\,]]\big],}$\\
 \indent $\mathsf{ Block[{\textcolor{blue}{\$MinPrecision}} = desirprec, {\textcolor{blue}{\$MaxExtraPrecision}} = \infty \}, gpcsetup[\textcolor{blue}{\emph{t}}\,]; asympt[\textcolor{blue}{\emph{t}}\,]+\sum_{j=1}^{j=lastn}storage[[\,j\,]] \,]\big]\Big]\Big]; }$\\
 \\
 \noindent \textcolor{teal}{\textsf{(*..........Using the code above, the short time algorithm is............*)}} \\
$\mathsf{gpcslow[\textcolor{blue}{\emph{t}}\,\_?NumericQ] := Quiet\big[If[\textcolor{blue}{\emph{t}}\, \leq \beta, 0, Block[\{{\textcolor{blue}{\$MinPrecision}} = 0, {\textcolor{blue}{\$MaxExtraPrecision}} = \infty\}, gpcshort[\textcolor{blue}{\emph{t}}\,]];}$\\
 \indent $\mathsf{Block[\{{\textcolor{blue}{\$MinPrecision}} = desirprec, {\textcolor{blue}{\$MaxExtraPrecision}} = \infty\}, gpcshort[\textcolor{blue}{\emph{t}}\,]]; \sum_{j=1}^{j=lastn}storage[[\,j\,]]\,]\big];}$\\
 \\
 \noindent \textcolor{teal}{\textsf{(*...........Using the code above, the long time algorithm is.............*)}} \\
$\mathsf{gpclong[\textcolor{blue}{\emph{t}}\,\_?NumericQ] := Quiet\Big[If[\textcolor{blue}{\emph{t}}\, \leq \beta, 0, }$\\
 \indent $\mathsf{Block\big[\{{\textcolor{blue}{\$MinPrecision}} = desirprec, {\textcolor{blue}{\$MaxExtraPrecision}} = \infty\}, gpcsetup[\textcolor{blue}{\emph{t}}\,]; asympt[\textcolor{blue}{\emph{t}}\,]+\sum_{j=1}^{j=lastn}storage[[\,j\,]]\big]\Big]\Big];}$\\
 \par}

\end{onehalfspacing}
\begin{small}
\bibliographystyle{myvancouver2}
\bibliography{Reff_GPC2}
\end{small}

\end{document}